\documentclass[11pt]{article}
\usepackage{geometry}                
\usepackage{graphicx}
\usepackage{amssymb}
\usepackage{latexsym}
\usepackage{amsopn}
\usepackage{amsthm,amsmath,amssymb} 
\usepackage{verbatim}
\usepackage{color}
\usepackage{epstopdf}
\newtheorem{Theorem}{Theorem}[section] 

\newtheorem{Lemma}[Theorem]{Lemma}

\newtheorem{Corollary}[Theorem]{Corollary} 

\theoremstyle{definition} 
\newtheorem{Definition}[Theorem]{Definition}

\newtheorem{Remark}[Theorem]{Remark}

\def\eps{\varepsilon}
\def\vu{{\bf u}}
\def\vv{{\bf v}}

\title{Double obstacle phase field approach to an inverse problem for a discontinuous diffusion coefficient}

\author{Klaus Deckelnick\footnote{Institut f\"ur Analysis und Numerik,
Otto-von-Guericke-Universit\"at Magdeburg, Universit\"atsplatz 2, 39106 Magdeburg, Germany.},
Charles M. Elliott\footnote{Mathematics Institute, University of Warwick, Coventry CV4 7AL, UK.}
and Vanessa Styles\footnote{Department of Mathematics,
University of Sussex, Brighton BN1 9RF, UK.} }

\date{}

\begin{document}
\maketitle

\begin{abstract}
  We propose a double obstacle phase field approach to the recovery of piece-wise constant diffusion coefficients
  for elliptic partial differential equations. The approach to this inverse problem is that of optimal control in which we have a 
quadratic  fidelity term to which we add a perimeter regularisation weighted by a parameter $\sigma$. This yields a 
functional which is optimised over a set of diffusion coefficients subject to a state equation which is the underlying 
elliptic PDE. In order to derive a problem which is amenable to computation 
  the perimeter functional is relaxed using a gradient energy functional together with an obstacle potential in 
  which there is an interface parameter $\epsilon$. This phase field  approach is justified by 
proving  $\Gamma-$convergence to the  functional  with perimeter regularisation as $\epsilon\rightarrow 0$. 
The computational approach is based on a finite 
element approximation. This discretisation is shown to converge in an appropriate way 
to the solution of the phase field problem. We derive an iterative method which 
is shown to yield an energy decreasing sequence 
 converging to a discrete critical point. The efficacy of the approach is illustrated with numerical
 experiments. 
\end{abstract}

\section{Introduction}
Many applications lead to mathematical models involving elliptic equations with piece-wise constant discontinuous coefficients. Frequently  
the interfaces  across which the coefficients  jump are completely unknown. A common approach  for  the 
identification of these coefficients is to make observations of the field variables solving the equations and use 
these values  in an attempt to determine the coefficients by formulating an inverse problem for the coefficients. 
This is generally ill posed and in applications it is usual to  use a  fidelity to the observations functional together with a regularisation of the coefficients. In this paper we use a regularisation of the  coefficients by employing the perimeter of the jump sets of the coefficients.

\subsection{Model problem}
To fix ideas we consider the following model elliptic problem:
\begin{eqnarray}\label{stateeq}
- \nabla \cdot \bigl(a \nabla y \bigr) & = & 0 \quad \mbox{ in } \Omega \\
a\frac{\partial y}{\partial \nu} & = & g \quad \mbox{ on } \partial \Omega, \label{statebv}
\end{eqnarray}
where $\Omega$ is a bounded domain in $\mathbb R^{d} \; (d=2,3)$, $g$ is given boundary data with zero 
mean
\begin{equation}  \label{zeromean}
\int_{\partial \Omega}g=0
\end{equation}
and $a$ is an isotropic  diffusion (conductivity) coefficient. 
We suppose that the diffusion coefficient takes one of the
$r$ positive values $a_1,\ldots, a_r$.  Our interest is in modelling a geometrical inverse problem 
concerning the determination of   the regions in which the material diffusion coefficient takes these values.  
Our problem then is to determine the sets $E_i=\lbrace x \in \Omega \, | \, 
a(x)=a_i \rbrace$ given observations of the solution $y$ of the elliptic boundary value problem (\ref{stateeq}), 
(\ref{statebv}). In the case of $r=2$, under constraints on the nature of the domains and boundary conditions,  uniqueness and stability results have been proved in \cite{BelFriIsa92, AleIsaPow95}. In this context see also  \cite{HetRun98}. 

A standard approach is to minimise a fidelity functional
\begin{displaymath}
J_{fid}(\mathcal E):= ||y_{\mathcal E}-y_{obs}||_{\mathcal O}^{2}
\end{displaymath}
over an appropriate  class of partitions  $\mathcal E =(E_i)_{i=1}^r$ of $\Omega$,  where $y_{\mathcal E}$ denotes the solution of the state or 
forward equation (\ref{stateeq}), (\ref{statebv}) with diffusion coefficient $a(x)=a_i, x \in E_i,i=1,\ldots,r$. Furthermore,
$\mathcal O$ is an appropriate space of observations and $y_{obs} \in \mathcal O$ is given.  In general this problem is  ill-posed and is typically regularised by adding a 
Tikhonov regularisation functional. A numerical approach without regularisation is proposed in \cite{HetRun98, KolArrLio99}.

\subsection{Geometric regularisation}
In this setting it has been considered appropriate to use perimeter regularisation, \cite{MUmSha89, ItoKunLi01}
\begin{displaymath}
J_{reg}(\mathcal E)=\hat{\sigma} \sum_{i=1}^r \mathcal H^{d-1}(\partial E_i \cap \Omega), \quad \mathcal E=(E_i)_{i=1}^r,
\end{displaymath}
where the regularisation parameter $\hat{\sigma}$ is positive. 
Minimisers of
\begin{displaymath}
J(\mathcal E):=J_{fid}(\mathcal E)+J_{reg}(\mathcal E)
\end{displaymath}
are then typically sought in the set of Caccioppoli partitions into $r$ components, i.e. partitions 
$\mathcal E=(E_i)_{i=1}^r$ of $\Omega$ with
$\mathcal H^d(E_i \cap E_j)=0, i \neq j, \; \mathcal H^d \bigl(\Omega \setminus \bigcup_{i=1}^r E_i \bigr)=0$ for which 
$u_i:=\chi_{E_i}$ belongs to $BV(\Omega), i=1,\ldots,r$. Thus, a Caccioppoli partition 
corresponds to a function $\vu=(u_1,\ldots,u_r) \in 
BV(\Omega,\lbrace e_1,\ldots,e_r \rbrace)$, where $e_1,\ldots,e_r$ are the unit vectors in $\mathbb{R}^r$. We can then
write the  regularisation functional in terms of $\vu$ as follows: 
\begin{displaymath}
J_{reg}(\vu)= \hat{\sigma} \sum_{i=1}^r \int_{\Omega} | D u_i |.
\end{displaymath}
Here,
$\int_{\Omega} | Du_i | $ is the total variation of the vector--valued Radon measure $Du_i$.
Before we rewrite the fidelity term
let us introduce the Gibbs simplex
\begin{displaymath}
\Sigma:= \lbrace {\bf y} \in \mathbb{R}^r \, | \, y_i \geq 0, i=1,\ldots,r, \; \sum_{i=1}^r y_i=1 \rbrace
\end{displaymath}
and observe that $e_1,\ldots,e_r$ are the corners of $\Sigma$. Consider the set 
\begin{displaymath}
X:= \lbrace \vu: \Omega \rightarrow \mathbb{R}^r \, | \, \vu \mbox{ is measurable and } \vu(x) \in \Sigma \mbox{ a.e. in } \Omega
\rbrace
\end{displaymath}
endowed with the $L^1$--norm and define for $\vu \in X$
\begin{equation} \label{diffcoeff}
a(\vu):= \sum_{i=1}^r a_i u_i
\end{equation}
and by $S(\vu)$ the solution of (\ref{stateeq}), (\ref{statebv}) with diffusion coefficient $a(\vu)$. 

We set $\hat{\sigma}= \frac{\pi}8\sigma$ for later convenience. The constant $\pi/8$  arises from the form of the phase field relaxation used in (\ref{eq:aq}), see (\ref{defF}).  

Problem ({\bf PGR}) is then to  seek
minimizers of the functional $J:X \rightarrow \mathbb{R} \cup \lbrace \infty \rbrace$ given by
\begin{displaymath}
J(\vu):= \left\{
\begin{array}{cl}
\displaystyle \frac{1}{2} ||S(\vu)-y_{obs}||^{2}_{\mathcal O} + \sigma \frac{\pi}{8} \sum_{i=1}^r \int_{\Omega} | Du_i | &, \mbox{ if } 
\vu \in BV(\Omega, \lbrace e_1,\ldots,e_r \rbrace) \cap X; \\[2mm]
\infty &, \mbox{ otherwise}.
\end{array}
\right.
\end{displaymath}

In this problem the fidelity term is  non-convex because of the nonlinearity of the state 
solution operator $S(\cdot)$ with respect to the coefficient $a(\vu)$. 
Also a feature of this natural geometric regularisation approach is that the regularisation functional is non-convex. 
This is reflected in the fact that $\vu$ only takes one of the
values $e_1,\ldots,e_r$ which leads to a non--convex constraint.

\subsection{Double obstacle phase field approach}
We shall  consider a suitable 
phase field approximation of the above regularisation  which involves gradient energies and functions that map into the Gibbs simplex.
 In this approximation we  relax the non-convex constraint $\vu(x) \in \lbrace e_1,\ldots,e_r \rbrace$ by introducing the set 
\begin{displaymath}
\mathcal K:= \lbrace \vu \in H^1(\Omega,\mathbb{R}^r) \, | \, \vu(x) \in \Sigma \mbox{ a.e. in } \Omega \rbrace
\end{displaymath}
and approximate $J$ by the sequence of functionals $J_{\epsilon}:X \rightarrow \mathbb{R} \cup \lbrace \infty \rbrace, \epsilon>0$ with
\begin{equation}
J_{\epsilon}(\vu):= \left\{
\begin{array}{cl}
\displaystyle \frac12 ||S(\vu)-y_{obs}||^{2}_{\mathcal O} + \sigma \int_\Omega \bigl( \frac{\eps}2| D \vu |^2 
+\frac1{2\eps}(1 - | \vu |^2) \bigr)dx &, \mbox{ if } \vu \in \mathcal K; \\
\infty &, \mbox{ otherwise}.
\end{array}
\right.
\label{eq:aq}
\end{equation}
Here, $\int_{\Omega} | D \vu |^2 dx = \sum_{i=1}^r \int_{\Omega} | \nabla u_i |^2 dx$ and 
for $\vu\in \mathcal K$ we have $ \int_\Omega (1 - | \vu |^2) =  \sum_{i=1}^r \int_{\Omega} u_i(1-u_i)$.  Problem ({\bf PDO}) is then to seek minimisers of $J_{\eps}$.
We refer to this approach as a double obstacle phase field model because of the constraints $0 \leq u_i \leq 1$ on the components of the phase field
vector $\vu$. The parameter $\epsilon$ is a measure of the thickness of a diffuse interface separating two sets on which the 
diffusion coefficient is constant.
The Cahn--Hilliard type energy 
\begin{displaymath}
\int_\Omega \bigl( \frac{\eps}2| \nabla u |^2 +\frac1{2\eps}(u - u^2) \bigr)dx
\end{displaymath}
is well established as an approximation of the perimeter functional, see e.g.
\cite{BloEll93-a,BloEll91,BelPaoVer90}. Note that the regularisation  remains   non-convex through the quadratic Cahn-Hilliard functional even though the constraint set is convex.
Let us remark that such a phase field model has recently been used in a binary recovery problem, see
\cite{BreDedEll13}.

\vspace{2mm}

{\color{black} Note that we view ({\bf PGO}) as having  just one regularisation parameter $\sigma$. 
The $\varepsilon$ parameter in ({\bf PGO}) may be viewed 
as a way of providing an approximation of ({\bf PGR}) which is computationally accessible. }


\subsection{Other approaches}

There have been attempts to solve the recovery problem without regularisation of the interfaces 
across which the diffusion coefficients jump. Formally one can write down variations of the fidelity functional with respect to variations of the interfaces. For example see \cite{HetRun98}.  In particular the interfaces can be associated with particular level sets of level set  functions
which have to be determined. 
We refer to \cite{San96,KolArrLio99, DorMilRap00, Bur03} for numerical implementations.
The use of level set descriptions of the interfaces in the context of perimeter regularisations is described in 
\cite{AmeBurHac04, Hac07, HegCanRic13}. Related to this is the use of total variation of a  regularised Heaviside function with argument being  a level set function, \cite{DeCLeiTai09,VdDAscLei10}. {\color{black} In \cite{multi-bang} the authors consider the distributed control of linear
elliptic systems in which the control variable should only take on a finite
number of values. To this purpose they introduce a combination of $L^2$ and
$L^0$--type penalties whose Fenchel conjugates allow the derivation of a
primal--dual optimality system with a unique solution. A suitable adaption
of this approach could be an alternative way to attack the inverse problem
considered in the present paper.}

In the  different context of image segmentation parametric description of curves have been used in 
conjunction with perimeter regularisation, \cite{BenGar14,TsaYezWil01}.

On the other hand  \cite{ChaTai03,TaiLi07, NieTaiAan07} use total variation  regularisation and relax the constraints that the indicator functions take just two values.

\subsection{Applications}
Our model problem is an example of the 
 identification of a  coefficient in an elliptic equation. This problem arises in many applications.  For example, 
 a fundamental issue in the use of mathematical models of flow in porous media is that the 
geological features which determine the permeability are unknown. In geology a facies is a body of rock 
with  specific characteristics. In our model problem $y$ is the pressure or hydraulic head associated with a 
fluid (for example, oil or water) occupying the reservoir or acquifer 
$\Omega$ and $a$ is the permeability of the rock. We assume  that the permeability is isotropic and is 
piece-wise constant. The domains $E_i=\lbrace x \in \Omega \, | \, 
a(x)=a_i \rbrace,~i=1,2,...,r$ model the decomposition of the reservoir $\Omega$ into facies whose location is 
unknown. The goal is to use observations of the pressure to determine the 
geometrical decomposition of the reservoir with respect to these facies, \cite{DorVil08, IglMcL11, IglLinStu14}.

Such problems also arise in imaging.
For example, electric impedance tomography, \cite{CheIsaNew99, DorMilRap00, BoyAdlLio12}, is the determination of the 
conductivity distribution in the interior of a 
domain using  observations of current and potential. Here $y$ is the electric potential and $a$ is a conductivity 
which takes different values in unknown interior domains. In medical imaging the shape and size of interior 
domains may be inferred from the variation of the conductivity.

\subsection{Outline and contributions of the paper}
\begin{itemize}
\item
In Section 2 we introduce the functionals $J_{\epsilon}$ and prove that 
they $\Gamma$--converge to $J$. Furthermore,
we show that $J_{\epsilon}$ has a minimum and derive a necessary first order condition. This establishes that problems ({\bf PGR}) and ({\bf PDO}) 
have solutions.
\item
The optimisation problem in Section 2 is infinite-dimensional. In order to carry out numerical calculations we employ a finite 
element spatial discretisation.  This is derived in Section 3 and we prove convergence results for absolute minimizers and critical
points as the mesh size tends to zero. This establishes that the inverse problems  ({\bf PGR}) and ({\bf PDO}) can be approximated by 
something computable.
\item
Section 4 is devoted to formulating  an iterative  scheme for finding critical points of the functional associated with the discrete 
optimisation problem. The method is based on a semi-implicit time discretisation of a parabolic variational inequality which is a 
gradient flow for the energy. In this finite dimensional setting we prove  a global convergence result for the iteration.
\item
Finally in Section 5 we illustrate the applicability of the method with some numerical examples.

\end{itemize}

\setcounter{equation}{0}

\section{Problem formulation}

\subsection{
State equation}
Let $\Omega \subset \mathbb{R}^d$ be a bounded domain with a Lipschitz boundary. We suppose that $g \in L^2(\partial \Omega)$ 
satisfying (\ref{zeromean}) and $y_{obs} \in \mathcal O$ are given functions. Here, $\bigl( \mathcal O,(\cdot,\cdot)_{\mathcal O} \bigr)$
is a Hilbert space with the property 
that $H^1(\Omega)$ is compactly embedded in $\mathcal O$. Furthermore we assume 
that  the following Poincar\'{e} inequality
\begin{equation} \label{poincare}
\displaystyle 
||\eta-\mathcal M_{\mathcal O}(\eta)||\le C_p||\nabla \eta||, \qquad \eta \in H^1(\Omega)
\end{equation}
holds, where $||\cdot||$ denotes the $L^{2}(\Omega)$ norm and  $\mathcal M_{\mathcal O}(\eta)$ denotes the mean value of $\eta$ with 
\begin{displaymath}
\mathcal M_{\mathcal O}(\eta):=(\eta,1)_{\mathcal O}/||1||^{2}_{\mathcal O}, \quad \eta \in \mathcal O.
\end{displaymath}
Typical examples are  $\mathcal O =L^{2}(\Omega)$ or $L^{2}(\partial \Omega)$ 
representing  either bulk measurements or boundary observations of the solution of the state equation. \\
For a given $\vu \in X$ we denote by $y=S(\vu) \in H^1(\Omega)$ the unique weak solution of the Neumann problem
\begin{eqnarray}
-\nabla\cdot (a(\vu)\nabla y)&=& 0 \quad \mbox{ in }\Omega \label{state1}\\
a(\vu)\frac{\partial y}{\partial \nu}& =& g \quad \mbox{ on }\partial \Omega \label{state2}
\end{eqnarray}
with  $\mathcal M_{\mathcal O}(y)=
\mathcal M_{\mathcal O}(y_{obs})$ 
in the sense that
\begin{equation}  \label{stateweak}
\displaystyle
\int_{\Omega} a(\vu) \nabla y \cdot \nabla \eta dx = \int_{\partial \Omega} g \eta do \qquad \forall \eta \in H^1(\Omega).
\end{equation}
Here, $a(\vu)$ is given by (\ref{diffcoeff}), where we note that
\begin{equation}  \label{a1a2}
a_{min} \leq a(\vu) \leq a_{max} \quad \mbox{ a.e. in } \Omega, 
\mbox{ uniformly in } \vu \in X,
\end{equation}
where $a_{min}:=\min(a_1,\ldots,a_r), \, a_{max}:= \max(a_1,\ldots,a_r)$.
Observe that $S$ is a nonlinear operator because of the bilinear relation between $a(\vu)$ and
$y$ in (\ref{stateweak}). Using (\ref{poincare}) together with the fact that $\mathcal M_{\mathcal O}(y)=
\mathcal M_{\mathcal O}(y_{obs})$ we infer that the solution $y=S(\vu)$ satisfies
\begin{displaymath}
\Vert y \Vert  \leq  \Vert y - \mathcal M_{\mathcal O}(y)\Vert + 
| \Omega |^{\frac{1}{2}} | \mathcal M_{\mathcal O}(y_{obs}) |  \leq C_p \Vert \nabla y \Vert+ 
\frac{ | \Omega |^{\frac{1}{2}}}{||1||_{\mathcal O}} \Vert y_{obs}\Vert_{\mathcal O}.
\end{displaymath}
If we combine this estimate with the choice $\eta=y$ in (\ref{stateweak}) and use (\ref{a1a2}) as well as the continuous embedding 
$H^1(\Omega) \hookrightarrow L^2(\partial \Omega)$ 
we deduce that
\begin{equation} \label{stateap}
\displaystyle 
\Vert S(\vu) \Vert_{H^1(\Omega)} \leq c(a_{min},\Omega) \bigl( \Vert g \Vert_{L^2(\partial \Omega)} + \Vert y_{obs} \Vert_{\mathcal O}\bigr) 
\quad \mbox{ uniformly in } \vu \in X.
\end{equation}
We see that the problem of observing $y$ given $\vu$ is well formulated because
 $$S:X \rightarrow \mathcal O~~\mbox{ is continuous}$$
  which is a consequence of the following lemma.
\begin{Lemma} \label{scont} $S:X \rightarrow H^1(\Omega)$ is continuous.
\end{Lemma}
{\it Proof.} Let $\vu \in X$ and $(\vu_k)_{k \in \mathbb{N}}$ a sequence in $X$ with $\vu_k \rightarrow \vu$ in
$L^1(\Omega,\mathbb{R}^r), k \rightarrow \infty$. Since $0 \leq u_{k,i} \leq 1, i=1,\ldots,r$ we may assume by passing to a subsequence if necessary that
$\vu_k \rightarrow \vu$ in $L^2(\Omega,\mathbb{R}^r)$ and a.e. in $\Omega$. Abbreviating $y=S(\vu),y_k=S(\vu_k)$ we have
for $\eta \in H^1(\Omega)$
\begin{displaymath}
\int_{\Omega} a(\vu_k) \nabla (y_k-y) \cdot \nabla \eta dx = \int_{\Omega} (a(\vu) - a(\vu_k)) \nabla y \cdot \nabla \eta dx.
\end{displaymath}
Choosing $\eta = y_k -y$ we deduce with the help of (\ref{a1a2}) and (\ref{diffcoeff})
\begin{displaymath}
a_{min} \Vert \nabla (y_k -y) \Vert \leq a_{max} \Bigl( \int_{\Omega} | \vu_k - \vu |^2 | \nabla y |^2 dx \Bigr)^{\frac{1}{2}}
\rightarrow 0, k \rightarrow \infty
\end{displaymath}
by the dominated convergence theorem because
\begin{displaymath}
|\vu_k - \vu |^2 | \nabla y |^2 \rightarrow 0 \mbox{ a.e. in } \Omega, \quad   | \vu_k - \vu |^2 | \nabla y |^2 
\leq r | \nabla y |^2 \mbox{ a.e in } \Omega \mbox{ and } | \nabla y |^2 \in L^1(\Omega).
\end{displaymath}
Since $\mathcal M_{\mathcal O}(y_k-y)=0$ we deduce
with the help of (\ref{poincare}) that $S(\vu_k)=y_k \rightarrow y=S(\vu)$ 
in $H^1(\Omega)$. \qed

\subsection{
$\Gamma$--convergence and existence of minimizers}

The use of $J_{\epsilon}$ in the minimization of $J$ is justified by the following $\Gamma$--convergence result.

\begin{Theorem} \label{gamma-conv}
The functionals $J_{\epsilon}$ $\Gamma$--converge to $J$ in $X$.
\end{Theorem}
{\it Proof.} Let us write $J_{\epsilon}(\vu)=G(\vu)+\sigma F_{\epsilon}(\vu)$, where 
$G(\vu)=\frac12 ||S(\vu)-y_{obs}||^{2}_{\mathcal O}$ is continuous as a consequence of Lemma \ref{scont} and the embedding of
$H^1(\Omega)$ into $\mathcal O$. 
In Theorem \ref{feps} in the Appendix we show that
\begin{equation} \label{defF}
\displaystyle
F_{\epsilon} \overset{\Gamma}{\rightarrow} F, \mbox{ where } F(\vu) = \left\{
\begin{array}{cl}
\displaystyle  \frac{\pi}{8} \sum_{i=1}^r \int_{\Omega} | Du_i | &, \mbox{ if } 
\vu \in BV(\Omega, \lbrace e_1,\ldots,e_r \rbrace) \cap X; \\[2mm]
\infty &, \mbox{ otherwise}.
\end{array}
\right.
\end{equation}
Using  Remark 1.7 in \cite{Bra02} we infer that $J_{\epsilon} \overset{\Gamma}{\rightarrow} G+\sigma F=J$. \qed

\begin{Theorem}  \label{existcont}
The minimization problem $\min_{\vv \in X} J_{\epsilon}(\vv)$ has a solution $\vu_{\epsilon} \in \mathcal K$.
\end{Theorem}
{\it Proof.} Let $(\vu_k)_{k \in \mathbb{N}} \subset \mathcal K$ be a minimizing sequence, $J_{\epsilon}(\vu_k) \searrow 
\inf_{\vv \in X} J_{\epsilon}(\vv)$. Since 
$(\vu_k)_{k \in \mathbb{N}}$ is bounded in $H^1(\Omega,\mathbb{R}^r)$ there exists a subsequence, again 
denoted by $(\vu_k)_{k \in \mathbb{N}}$,
and $\vu_{\epsilon} \in H^1(\Omega,\mathbb{R}^r)$ such that
\begin{displaymath}
\vu_k \rightharpoonup \vu_{\epsilon} \mbox{ in } H^1(\Omega,\mathbb{R}^r), \quad \vu_k \rightarrow \vu_{\epsilon} \mbox{ in } 
L^2(\Omega,\mathbb{R}^r) \mbox{ and a.e. in } \Omega.
\end{displaymath}
In particular, $\vu_{\epsilon} \in \mathcal K$. Lemma \ref{scont} implies that $S(\vu_k) \rightarrow S(\vu_{\epsilon})$ 
in $\mathcal O$ which combined with the weak lower semicontinuity
of the $H^1$-seminorm shows that $\vu_{\epsilon}$ is a minimum of $J_{\epsilon}$. \qed

\begin{Corollary}
Let $(u_{\epsilon})_{\epsilon>0}$ be a sequence of minimizers of $J_{\epsilon}$. Then there exists a sequence $\epsilon_k \rightarrow
0, k \rightarrow \infty$ and $\vu \in BV(\Omega;\lbrace e_1,\ldots,e_r \rbrace) \cap X$ such that $\vu_{\epsilon_k}
\rightarrow \vu$ in $L^1(\Omega,\mathbb{R}^r)$ and $\vu$ is a minimum of $J$.
\end{Corollary}
{\it Proof.} By Corollary \ref{compact} in the Appendix there exists a sequence $\epsilon_k \rightarrow
0, k \rightarrow \infty$ and $\vu \in BV(\Omega;\lbrace e_1,\ldots,e_r \rbrace) \cap X$ such that $\vu_{\epsilon_k}
\rightarrow \vu$ in $L^1(\Omega,\mathbb{R}^r)$. It is well--known that the $\Gamma$--convergence of $J_{\epsilon_k}$ to $J$
implies that $\vu$ is a minimum of $J$. \qed

\vspace{3mm}

 \subsection{Necessary first order condition for the phase field recovery}
In order to derive the necessary first order conditions for a minimum of $J_{\eps}$ we consider $\mathcal K$ as a subset
of $L^{\infty}(\Omega,\mathbb{R}^r)$. Similarly as in \cite{BlaFarGar13}, Section 3, one can prove that the solution operator
$S:L^{\infty}(\Omega,\mathbb{R}^r) \supset \mathcal K \rightarrow H^1(\Omega)$ is Fr\'{e}chet differentiable with
$\tilde{y}=S'(\vu){\bf w}, {\bf w} \in L^{\infty}(\Omega,\mathbb{R}^r)$ being given as the solution of
\begin{equation} \label{eq:1}
 \int_\Omega a(\vu)\nabla \tilde{y} \cdot \nabla \eta dx = -\int_\Omega  a({\bf w}) \nabla S(\vu) \cdot \nabla\eta dx \qquad
\forall\eta \in H^1(\Omega)
\end{equation}
with $\mathcal M_{\mathcal O}(\tilde{y}) =0$. As a result, $J_{\eps}$ is Fr\'{e}chet differentiable on
$\mathcal K \subset L^{\infty}(\Omega,\mathbb{R}^r) \cap H^1(\Omega,\mathbb{R}^r)$ with 
\begin{equation} \label{eq:2}
\displaystyle
J_{\eps}'(\vu){\bf w} = \bigl(S(\vu)-y_{obs},S'(\vu) {\bf w} \bigr)_{\mathcal O}   
+\sigma  \int_{\Omega}\bigl(  \eps D \vu \cdot D {\bf w} -
 \frac1{\eps} \vu \cdot {\bf w} \bigr) dx
\end{equation}
for ${\bf w} \in L^{\infty}(\Omega,\mathbb{R}^r) \cap H^1(\Omega,\mathbb{R}^r)$. 
In order to avoid the evaluation of $S'(\vu){\bf w}$ in (\ref{eq:2}) we work as usual with a dual problem: 
Find $p \in H^1(\Omega)$ such that $\mathcal M_{\mathcal O}(p)=0$ 
and
\begin{equation}
\int_\Omega a(\vu)\nabla p \cdot \nabla \eta dx = \bigl( S(\vu)-y_{obs}, \eta \bigr)_{\mathcal O} ~~\forall\eta \in H^1(\Omega),
\label{eq:dual}
\end{equation}
where we note that the solvability condition $\bigl( S(\vu)-y_{obs}, 1 \bigr)_{\mathcal O}=0$ is satisfied. As a result we obtain
from (\ref{eq:2}), (\ref{eq:dual}) and (\ref{eq:1}) 
\begin{eqnarray}
J_{\eps}'(\vu){\bf w}  
&   =  & \int_{\Omega}a(\vu)\nabla p \cdot \nabla [S'(\vu){\bf w]} dx 
+\sigma \int_{\Omega} \bigl( \eps D \vu \cdot D {\bf w} - \frac1{\eps} \vu \cdot {\bf w} \bigr) dx \nonumber
\label{eq:3} \\
& = & - \int_{\Omega} a({\bf w}) \, \nabla S(\vu) \cdot \nabla p dx 
+\sigma \int_{\Omega} \bigl( \eps D \vu  \cdot D {\bf w} - \frac1{\eps} \vu \cdot {\bf w} \bigr) dx. \nonumber
\end{eqnarray}
At a minimum $\vu$ of $J_{\epsilon}$ we have $J_{\eps}'(\vu)(\vv-\vu) \geq 0$ for all $\vv \in \mathcal K$. Since 
$a(\vv-\vu)=a(\vv) - a(\vu)$ we therefore define: 

\begin{Definition}({\bf Phase field critical point})
Find $\vu\in \mathcal K$ such that for all $\vv \in \mathcal K$
\begin{equation} \label{eq:100}
\displaystyle 
\sigma \int_{\Omega} \bigl( \eps D \vu  \cdot D ( \vv - \vu)  - \frac1{\eps} \vu \cdot  (\vv - \vu) \bigr) dx
 -  \int_{\Omega} ( a(\vv) - a(\vu)) \, \nabla S(\vu) \cdot \nabla p dx  \geq 0.
\end{equation}
\end{Definition}

\begin{Remark} \label{parobsprob}
A natural  strategy to construct solutions of (\ref{eq:100}) and hence to find candidates 
for at least a local minimum of $J_{\eps}$ is to consider the following parabolic obstacle problem: Find $\vu(\cdot,t) \in \mathcal K, t \geq 0$
such that $\vu(\cdot,0)=\vu_0$ and 
\begin{displaymath}
(\vu_t, \vv-\vu)
+\sigma  \int_{\Omega} \bigl( \eps D \vu  \cdot D ( \vv - \vu) - \frac1{\eps} \vu \cdot (\vv - \vu) \bigr)dx 
- \int_{\Omega} (a(\vv) - a(\vu)) \, \nabla S(\vu) \cdot \nabla p dx \geq 0 
\end{displaymath}
for all $\vv\in \mathcal K$ and all $t>0$. Here, $p$ is the solution of (\ref{eq:dual}) for $\vu(\cdot,t)$ and
$\vu_0 \in \mathcal K$ is a suitably chosen initial function.
\end{Remark}
Inserting ${\bf v}=\vu(\cdot,t-\Delta t)$ into the above relation, dividing by $\Delta t$ and sending $\Delta t \rightarrow 0$
we formally find that
\begin{displaymath}
\Vert \vu_t \Vert^2 + J_{\epsilon}'(\vu)\vu_t \leq 0,
\end{displaymath}
so that $\displaystyle \frac{d}{dt} J_{\epsilon}(\vu(\cdot,t)) \leq 0$ and the value of the objective funtional
decreases during the evolution. If $\lim_{t \rightarrow \infty} u(\cdot,t)=:u_{\infty}$ exists, we expect $u_{\infty}$ to be
 a solution of (\ref{eq:100}).

\setcounter{equation}{0}
\section{Finite element approximation}
In what follows we assume that $\Omega$ is a polygonal (d=2) or  polyhedral (d=3) domain. Let us denote by $(\mathcal T_h)_{0<h \leq h_0}$ 
a regular triangulation of $\Omega$ and set
\begin{displaymath}
V_h = \lbrace \chi \in C^0(\bar{\Omega}) \,  | \, \chi_{|T} \in P_1(T) \mbox{ for all } T \in \mathcal T_h \rbrace\subset H^{1}(\Omega)
\end{displaymath}
as well as 
\begin{displaymath}
\mathcal K_h:= \lbrace {\bf \chi} \in (V_h)^r \, | \,  {\bf \chi}(x)  \in \Sigma, x \in \bar{\Omega} \rbrace \subset \mathcal K.
\end{displaymath}
Using the construction of the Cl\'{e}ment interpolation operator (\cite{Cle75}) it is not difficult to see that for every $\vu \in \mathcal K$
there exists a sequence $(\hat{\vu}_h)_{0 < h \leq h_0}$ with $\hat{\vu}_h \in \mathcal K_h$ such that 
\begin{equation} \label{approx}
\displaystyle \hat{\vu}_h \rightarrow \vu \mbox{ in } H^1(\Omega,\mathbb{R}^r) \mbox{ as } h \rightarrow 0.
\end{equation}
Furthermore, let $(y_{obs}^{h})_{0< h \leq h_0}$ be a sequence of functions $y_{obs}^{h}\in \mathcal O$ 
such that
\begin{equation}  \label{fh}
\displaystyle y^{h}_{obs} \rightarrow y_{obs} \mbox{ in }\mathcal O \mbox{ as } h \rightarrow 0. 
\end{equation}
For $\vu_h \in \mathcal K_h$ we denote by $y_h=S_h(\vu_h) \in V_h$ the solution of
\begin{equation}
\int_{\Omega} a(\vu_h)\nabla y_h \cdot  \nabla \chi dx  = \int_{\partial \Omega} g_h \chi do \quad \forall\chi \in V_h
\label{eq:stated}
\end{equation}
with $\mathcal M_{\mathcal O}(y_h) = \mathcal M_{\mathcal O}(y^h_{obs})$. Here $g_h: \partial \Omega \rightarrow \mathbb{R}$
is a piecewise linear, continuous approximation to  $g$ satisfying 
\begin{equation}  \label{gh}
\displaystyle
\int_{\partial \Omega} g_h do=0 \mbox{ and } g_h \rightarrow g \mbox{ in } L^2(\partial \Omega) \mbox{ as } h \rightarrow 0.
\end{equation}
In the same way as in (\ref{stateap}) one can prove that
\begin{equation}  \label{eststated}
\displaystyle \Vert S_h(\vu_h) \Vert_{H^1} \leq c \bigl( \Vert g_h \Vert_{L^2(\partial \Omega)} + \Vert y^h_{obs} \Vert_{\mathcal O}\bigr)  
\leq c
\quad \mbox{ uniformly in } \vu_h \in \mathcal K_h,
\end{equation}
where the constant $c$ is independent of $h$ in view of (\ref{fh}) and (\ref{gh}).

\begin{Lemma} \label{shcont} Let $(h_k)_{k \in \mathbb{N}}$ be a sequence with $\lim_{k \rightarrow \infty} h_k =0$ and
$\vu_{h_k} \in \mathcal K_{h_k}$ with
$\vu_{h_k} \rightarrow \vu$ in $L^1(\Omega,\mathbb{R}^r)$. Then $S_{h_k}(\vu_{h_k}) \rightarrow S(\vu)$ in $H^1(\Omega), k 
\rightarrow \infty$.
\end{Lemma}
{\it Proof.} Let $\vu_k=\vu_{h_k}, y_k=S_{h_k}(\vu_k)$ and $ y=S(\vu)$. By passing to a subsequence if necessary we may assume 
in addition that
$\vu_k \rightarrow \vu$ a.e. in $\Omega$. Choose a sequence $\hat{y}_k \in V_{h_k}$ 
such that $\hat{y}_k \rightarrow  y$ in $H^1(\Omega)$. Using (\ref{poincare}) we deduce
\begin{eqnarray}
\Vert y_k - \hat{y}_k \Vert_{H^1} & \leq & \Vert y_k - \hat{y}_k - \mathcal M_{\mathcal O}(y_k - \hat{y}_k) \Vert + | \Omega |^{\frac{1}{2}}
| \mathcal M_{\mathcal O}(y^{h_k}_{obs} - \hat{y}_k) | +  \Vert \nabla ( y_k - \hat{y}_k) \Vert \nonumber  \\
& \leq & c \Vert \nabla ( y_k - \hat{y}_k) \Vert + | \Omega |^{\frac{1}{2}} \bigl( | \mathcal M_{\mathcal O}(y^{h_k}_{obs} - y_{obs}) |
+ | \mathcal M_{\mathcal O}(y - \hat{y}_k) | \bigr) \nonumber \\
& \leq & c \Vert \nabla ( y_k - \hat{y}_k) \Vert + c ( \Vert y^{h_k}_{obs} - y_{obs} \Vert_{\mathcal O} + 
\Vert y - \hat{y}_k \Vert_{H^1}). \label{l2dif}
\end{eqnarray}
In order to estimate the first term we write
\begin{eqnarray*}
\lefteqn{ \int_{\Omega} a(\vu_k) \nabla (y_k - \hat{y}_k) \cdot \nabla \chi dx } \\
& = & \int_{\Omega} a(\vu_k) \nabla ( y - \hat{y}_k) \cdot \nabla \chi dx + \int_{\Omega} \bigl( a(\vu)-a(\vu_k) \bigr) \nabla y \cdot \nabla \chi dx 
+ \int_{\partial \Omega} (g_{h_k} - g) \chi do
\end{eqnarray*}
for all $\chi \in V_{h_k}$.  If we let $\chi=y_k-\hat{y}_k$ and take into account (\ref{l2dif}) we obtain
\begin{eqnarray*}
\Vert  y_k-\hat{y}_k \Vert_{H^1} &  \leq &  c \Vert y -\hat{y}_k \Vert_{H^1} + c \Bigl( \int_{\Omega} | \vu_k - \vu |^2 | \nabla y |^2 dx \Bigr)^{\frac{1}{2}} \\
& &  + c \bigl( \Vert g_{h_k} - g \Vert_{L^2(\partial \Omega)} +  \Vert y^{h_k}_{obs} - y_{obs} \Vert_{\mathcal O}\bigr) \rightarrow 0,
k \rightarrow \infty
\end{eqnarray*}
by (\ref{gh}) and (\ref{fh}). Here, the second integral is shown to converge to zero in the same way as in the proof of Lemma \ref{scont}. 
In conclusion,
$S_{h_k} (\vu_{h_k}) = (y_k-\hat{y}_k) + \hat{y}_k \rightarrow y=S(\vu)$ in $H^1(\Omega)$ and by a standard argument the whole sequence
converges. \qed

\vspace{3mm}

Using $S_h$ we define the following approximation $J_{\epsilon,h}: \mathcal K_h \rightarrow 
\mathbb{R}$ of $J_{\epsilon}$:
\begin{equation}
J_{\epsilon,h}(\vu_h):= \frac12||S_{h}(\vu_{h})-y^{h}_{obs}||^2_{\mathcal O}+ \sigma 
 \int_\Omega \bigl( \frac{\eps}2| D \vu_h|^2 
+\frac1{2\eps} (1 - | \vu_h |^2) \bigr)dx.
\label{eq:res}
\end{equation}

\begin{Theorem} \label{existdisc}
There exists $\vu_h \in \mathcal K_h$ such that $J_{\epsilon,h}(\vu_h) = \min_{\vv_h \in \mathcal K_h} 
J_{\epsilon,h}(\vv_h)$. Every sequence
$(\vu_{h_k})_{k \in \mathbb{N}}$ with $\lim_{k \rightarrow \infty} h_k=0$ has a subsequence that converges strongly in 
$H^1(\Omega,\mathbb{R}^r)$ and a.e. in $\Omega$ to a minimum of $J_{\epsilon}$.
\end{Theorem}
{\it Proof.} Since $X_h$ is finite-dimensional, the existence of a minimum of $J_{\epsilon,h}$ is straightforward. Next, let 
$\vu_k \in \mathcal K_{h_k}$ be a sequence with $\lim_{k \rightarrow \infty}h_k =0$ and 
$J_{\epsilon,h_k}(\vu_k)=\min_{\vv_h \in \mathcal K_{h_k}}J_{\epsilon,h_k}(\vv_h)$.
Since $(\vu_k)_{k \in \mathbb{N}}$ is bounded in $H^1(\Omega,\mathbb{R}^r)$, there exists a subsequence, again
denoted by $(\vu_k)_{k \in \mathbb{N}}$, and $\vu \in \mathcal K$ such that
\begin{equation} \label{conv1}
\displaystyle 
\vu_k  \rightharpoonup  \vu \mbox{ in }  H^1(\Omega,\mathbb{R}^r), \quad
\vu_k \rightarrow  \vu \mbox{ in } L^1(\Omega,\mathbb{R}^r) \mbox{ and a.e. in } \Omega.
\end{equation}
Furthermore, Lemma \ref{shcont} implies that
\begin{equation}  \label{conv2}
\displaystyle S_{h_k}(\vu_k) \rightarrow S(\vu) \mbox{ in } H^1(\Omega).
\end{equation}
We claim that $\vu$ is a minimum of $J_{\epsilon}$. To see this,  let $\vv \in \mathcal K$ be arbitrary and 
$\hat{\vv}_k \in \mathcal K_{h_k}$ a sequence
with $\hat{\vv}_k \rightarrow \vv$ in $H^1(\Omega,\mathbb{R}^r)$, see (\ref{approx}). Since
$J_{\epsilon,h_k}(\vu_k) \leq J_{\epsilon,h_k}(\hat{\vv}_k)$ we deduce from (\ref{conv1}), (\ref{conv2}) and again
Lemma \ref{shcont} that
\begin{displaymath}
J_{\epsilon}(\vu) \leq \liminf_{k \rightarrow \infty} J_{\epsilon,h_k}(\vu_k) \leq \limsup_{k \rightarrow \infty}
J_{\epsilon,h_k}(\vu_k) \leq \lim_{k \rightarrow \infty} J_{\epsilon,h_k}(\hat{\vv}_k) = J_{\epsilon}(\vv),
\end{displaymath}
so that $J_{\epsilon}(\vu) = \min_{\vv \in X} J_{\epsilon}(\vv)$. Furthermore, by repeating the above argument with
a sequence $\hat{\vu}_k \in \mathcal K_{h_k}$ such that
$\hat{\vu}_k \rightarrow \vu$ in $H^1(\Omega,\mathbb{R}^r)$ we infer in addition that
\begin{equation} \label{convobj}
\displaystyle \lim_{k \rightarrow \infty} J_{\epsilon,h_k}(\vu_k) = J_{\epsilon}(\vu).
\end{equation}
We use this relation to show that  $ \Vert D \vu_k \Vert^2 \rightarrow \Vert D \vu \Vert^2$. Namely, let us write
\begin{eqnarray*}
\frac{\sigma \epsilon}{2} \int_{\Omega} | D \vu_k|^2 dx & = &  J_{\epsilon,h_k}(\vu_k) - \frac{\sigma}{2 \epsilon} 
\int_{\Omega} (1 - | \vu_k |^2) dx  - \frac{1}{2} \Vert S_{h_k}(\vu_k) - y^{h_k}_{obs} \Vert^2_{\mathcal O}  \\ 
& \rightarrow &  J_{\epsilon}(\vu)  - \frac{\sigma}{2 \epsilon} \int_{\Omega} ( 1 - | \vu |^2 ) dx - \frac{1}{2} 
\Vert S(\vu) - y_{obs} \Vert^2_{\mathcal O} = \frac{\sigma \epsilon}{2} \int_{\Omega} | D \vu |^2 dx
\end{eqnarray*}
in view of (\ref{convobj}), (\ref{conv1}), (\ref{conv2}) and (\ref{fh}). Hence $\vu_k \rightarrow \vu$ in $H^1(\Omega,\mathbb{R}^r)$ and the
theorem is proved. \qed

\vspace{3mm}

In practice, rather than trying to locate a global minimum of $J_{\epsilon,h}$ one looks for admissible points $\vu_h$ that satisfy the
necessary first order condition 
\begin{equation} \label{discvar}
\displaystyle
J_{\epsilon,h}'(\vu_h)(\vv_h-\vu_h) \geq 0 \quad \mbox{  for all } \vv_h \in \mathcal K_h.
\end{equation}
A calculation analogous to (\ref{eq:3}) leads us to the following variational inequality: 
\begin{equation} \label{discvar1}
\displaystyle  \sigma  \int_{\Omega} \bigl( \eps D \vu_h \cdot D ( \vv_h - \vu_h) - \frac1{\eps} \vu_h \cdot (\vv_h - \vu_h) \bigr) dx
- \int_{\Omega} (a(\vv_h) - a(\vu_h)) \, \nabla y_h \cdot \nabla p_h dx  \geq 0
\end{equation}
for all $\vv_h \in \mathcal K_h$, where $y_h=S_h(\vu_h)$ and $p_h \in V_h$ with 
$\mathcal M_{\mathcal O}(p_{h})=0$
is the solution of 
the discrete adjoint problem: 
\begin{equation} \label{eq:duald}
\displaystyle
\int_\Omega a(\vu_h)\nabla p_h \cdot \nabla \chi dx = (y_{h}-y^{h}_{obs},\chi)_{\mathcal O}  \quad \forall\chi \in V_h.
\end{equation}

\begin{Theorem}
Let  $(\vu_{h_k})_{k \in \mathbb{N}}$ be a sequence of solutions of (\ref{discvar1}) with $\lim_{k \rightarrow \infty} h_k=0$. Then there
exists   a subsequence that converges strongly in $H^1(\Omega,\mathbb{R}^r)$ and
a.e. in $\Omega$ to a solution $\vu$ of (\ref{eq:100}).
\end{Theorem}
{\it Proof.} Let us abbreviate $\vu_k=\vu_{h_k}, y_k=S_{h_k}(\vu_k)$ and denote by $p_k \in V_{h_k}$ the solution of (\ref{eq:duald})
with $\vu_h=\vu_k$ and $y_h=y_k$. Using  (\ref{eststated}) and testing (\ref{eq:duald}) with $\chi = p_k$ we infer that 
\begin{displaymath}
\Vert y_k \Vert_{H^1} + \Vert p_k \Vert_{H^1}  \leq c \qquad \mbox{ uniformly in } k \in \mathbb{N}.
\end{displaymath}
Next, inserting $\vv_h \equiv  \frac{1}{r} \sum_{j=1}^r e_j$ into (\ref{discvar1}) we deduce
\begin{eqnarray*}
\sigma \eps  \int_{\Omega} | D \vu_k |^2 dx & \leq & \frac{\sigma}{\eps}   \int_{\Omega} | \vu_k |^2 dx  
+    \int_{\Omega}( a(\vu_k) - \frac{1}{r} \sum_{i=1}^r a_i) \,  \nabla y_k \cdot \nabla p_k dx  \\
& \leq & \frac{\sigma r}{\eps} | \Omega | + c \Vert \nabla  y_k \Vert \,  \Vert \nabla p_k \Vert \leq c.
\end{eqnarray*} 
Hence, there exists a subsequence, again
denoted by $(\vu_k)_{k \in \mathbb{N}}$, and $\vu \in \mathcal K$ such that
\begin{equation}  \label{conv3}
\displaystyle 
\vu_k \rightharpoonup \vu \mbox{ in } H^1(\Omega,\mathbb{R}^r), \quad \vu_k \rightarrow \vu \mbox{ in } 
L^1(\Omega,\mathbb{R}^r) \mbox{ and a.e. in } \Omega.
\end{equation}
Lemma \ref{shcont} implies that
\begin{equation}  \label{conv4}
\displaystyle y_k=S_{h_k}(\vu_k) \rightarrow S(\vu)=:y \mbox{ in } H^1(\Omega).
\end{equation}
Let $p \in H^1(\Omega), {\mathcal M}_{\mathcal O}(p)=0$ be the solution of (\ref{eq:dual}).
Choose $\hat{p}_k \in V_{h_k}$ with ${\mathcal M}_{\mathcal O}(\hat{p}_k)=0$ 
such that $\hat{p}_k \rightarrow  p$ in $H^1(\Omega)$ and write
\begin{eqnarray*}
\lefteqn{ 
\int_{\Omega} a(\vu_k) \nabla (p_k - \hat{p}_k) \cdot \nabla \chi dx = \int_{\Omega} a(\vu_k) \nabla ( p - \hat{p}_k) \cdot \nabla \chi dx } \\
& & + \int_{\Omega} \bigl( a(\vu)-a(\vu_k) \bigr) \nabla p \cdot \nabla \chi dx  
 + (y_k - y, \chi)_{\mathcal O} -  (y_{obs}^{h_{k}}-y_{obs},\chi)_{\mathcal O}
\end{eqnarray*}
for all $\chi \in V_{h_k}$. By choosing $\chi = p_k - \hat{p}_k$ and using (\ref{poincare}), (\ref{conv4}) and (\ref{fh}) we deduce 
\begin{displaymath}
\Vert p_k - \hat{p}_k \Vert_{H^1} \leq c \Vert \hat{p}_k - p \Vert_{H^1} + c \Bigl( \int_{\Omega} | \vu_k - \vu |^2 | \nabla p |^2 dx \Bigr)^{\frac{1}{2}}
+c \bigl( \Vert y_k - y \Vert_{\mathcal O} + \Vert y^{h_k}_{obs}- y_{obs} \Vert_{\mathcal O} \bigr)
\rightarrow 0
\end{displaymath}
which implies that $p_k \rightarrow p$ in $H^1(\Omega)$. \\
Let us next show that $\vu$ satisfies (\ref{eq:100}). Given $\vv \in \mathcal K$ there exists a sequence
$\hat{\vv}_k \in  \mathcal K_{h_k}$ such that $\hat{\vv}_k \rightarrow \vv$ in 
$H^1(\Omega,\mathbb{R}^r)$ and a.e. in $\Omega$. Then we have from (\ref{discvar1})
\begin{equation} \label{discvar2}
\displaystyle
\sigma  \int_{\Omega} \bigl( \eps D \vu_k \cdot D (\hat{\vv}_k - \vu_k) -\frac{1}{\eps} \vu_k \cdot (\hat{\vv}_k - \vu_k) \bigr) dx
- \int_{\Omega} (a(\hat{\vv}_k)-a(\vu_k)) \, \nabla y_k \cdot \nabla p_k dx  \geq 0.
\end{equation}
In order to examine the second term we write 
\begin{eqnarray}
\lefteqn{ \hspace{-1cm}
 \int_{\Omega}(a(\hat{\vv}_k) - a(\vu_k)) \, \nabla y_k \cdot \nabla p_k dx  - \int_{\Omega}(a(\vv) - a(\vu)) \, \nabla y \cdot \nabla p dx  } \label{conv5}  \\
& =  &   \int_{\Omega}(a(\hat{\vv}_k) - a(\vu_k))[ \nabla (y_k-y) \cdot \nabla p_k + \nabla y \cdot \nabla (p_k-p)]dx \nonumber \\
& & + \int_{\Omega}\bigl( (a(\hat{\vv}_k) - a(\vv)) -(a(\vu_k) - a(\vu)) \bigr) \, \nabla y \cdot \nabla p dx  
 \rightarrow  0, k \rightarrow \infty \nonumber
\end{eqnarray}
since $y_k \rightarrow y, p_k \rightarrow p$ in $H^1(\Omega)$ where we used again the dominated convergence theorem
for the second integral. By passing to the limit in (\ref{discvar2}) and
observing that $\int_{\Omega}  | D \vu |^2 dx \leq \liminf_{k \rightarrow \infty} \int_{\Omega} | D \vu_k |^2 dx$
we infer that $\vu$ satisfies ({\ref{eq:100}).  \\
Let us finally show that $\vu_k \rightarrow \vu$ in $H^1(\Omega,\mathbb{R}^r)$. Choose a sequence 
$\hat{\vu}_k \in \mathcal K_{h_k}$ such that
$\hat{\vu}_k \rightarrow \vu$ in $H^1(\Omega,\mathbb{R}^r)$. Inserting
$\vv_{h_k}=\hat{\vu}_k$ into (\ref{discvar1}) we obtain
\begin{displaymath}
\sigma \epsilon \int_{\Omega} | D \vu_k  |^2 dx 
 \leq  \sigma \epsilon \int_{\Omega} D \vu_k \cdot D \hat{\vu}_k  dx - \frac{\sigma}{\epsilon} \int_{\Omega}
 \vu_k \cdot ( \hat{\vu}_k - \vu_k)  dx  - \int_{\Omega}
\bigl( a(\hat{\vu}_k) - a(\vu_k) \bigr) \, \nabla y_k \cdot \nabla p_k dx 
\end{displaymath}
so that  (\ref{conv3}) and (\ref{conv5}) with $\hat{\vv}_k= \hat{\vu}_k$ imply that 
\begin{displaymath}
\limsup_{k \rightarrow \infty} \int_{\Omega} | D \vu_k  |^2 dx \leq \int_{\Omega} | D \vu |^2 dx.
\end{displaymath}
Hence $\int_{\Omega} | D \vu_k  |^2 dx \rightarrow  \int_{\Omega} | D \vu |^2 dx$, so that $D \vu_k \rightarrow D \vu$ in
$L^2$.  \qed

\setcounter{equation}{0}
\section{
An iterative scheme}
\subsection{Iterative method}
Let us consider  the following iteration, which can be seen as a time discretization of the parabolic obstacle 
problem introduced in Remark \ref{parobsprob}. Given
$\vu^n_h \in \mathcal K_h$ let $\vu^{n+1}_h \in \mathcal K_h$ be the solution of the problem 
\begin{eqnarray}
\lefteqn{  \int_{\Omega} (\vu_h^{n+1}-\vu^n_h) \cdot ( \vv_h - \vu_h^{n+1}) dx  - \tau_n  
 \int_{\Omega} (a(\vv_h) - a(\vu^{n+1}_h)) \, \nabla y_h^n \cdot  \nabla p_h^n dx }  \label{eq:flowd}  \\
& &  +\tau_n \sigma  \int_{\Omega} \bigl( \eps D \vu_h^{n+1} \cdot  D (\vv_h - \vu_h^{n+1}) -  
\frac1{\eps} \vu_h^{n} \cdot (\vv_h - \vu_h^{n+1}) \bigr) dx
\geq 0 \quad \forall \vv_h \in \mathcal K_h,  \nonumber
\end{eqnarray}
where $\tau_n>0$, $y^n_h= S_h(\vu^n_h)$ and  $p^n_h \in V_h$ solves the discrete dual problem
\begin{equation}
\int_{\Omega} a(\vu_h^n)\nabla p_h^n \cdot  \nabla \chi dx =(y_{h}^{n}-y^h_{obs},\chi)_{\mathcal O} ~~\forall\chi \in V_h \mbox{ with }
{\mathcal M}_{\mathcal O}(p_h^n)=0.  
\end{equation}
Note that $\vu^{n+1}_h$ is the unique solution of the convex minimization problem
\begin{displaymath}
\min_{\vv_h \in \mathcal K_h} \Bigl( \frac{1}{2} \Vert \vv_h - \vu^n_h \Vert^2 - \tau_n \int_{\Omega} a(\vv_h) \, 
\nabla y^n_h \cdot \nabla p^n_h dx
+ \tau_n \sigma \int_{\Omega} \bigl( \frac{\epsilon}{2} | D \vv_h |^2 dx  - \frac{1}{\epsilon} \vu^n_h \cdot \vv_h \bigr) dx \Bigr).
\end{displaymath} 
\subsection{
Convergence of the iterative method}
The following result shows that the objective functional decreases in the iteration provided the time steps $\tau_n$ satisfy
a suitable condition. In order to formulate it we define
\begin{equation} \label{poincare1}
\displaystyle \hat{a}:= \bigl( \sum_{i=1}^r a_i^2 \bigr)^{\frac{1}{2}}, \quad
\hat{c}:= \inf\Big\{ \,  \frac{ \int_{\Omega} | \nabla \eta |^2 dx}{\Vert \eta \Vert_{\mathcal O}^2} \, | \,
\eta \in H^1(\Omega) \setminus \lbrace 0 \rbrace, \mathcal M_{\mathcal O}(\eta)=0 \Big\}.
\end{equation}
Note that $\hat{c}\geq C_p^2>0$ in view of (\ref{poincare}).

\begin{Lemma} \label{monotonicity}
The sequence $(\vu^n_h)_{n \in \mathbb{N}_0}$  satisfies
\begin{displaymath}
\Vert \vu^{n+1}_h - \vu^n_h \Vert^2 + J_{\epsilon,h}(\vu^{n+1}_h) \leq J_{\epsilon,h}(\vu^n_h), \quad n \in \mathbb{N}_0,
\end{displaymath}
provided that 
\begin{equation}  \label{timestepcond}
\displaystyle \tau_n \leq \Bigl(1+ 
\frac{\hat{a}^2}{a_{min}} \Vert \nabla y^n_h \Vert_{L^{\infty}} \Vert \nabla p^n_h \Vert_{L^{\infty}}
+ \frac{\hat{a}^2}{a_{min}^2}\frac{1}{2 \hat{c}} \Vert \nabla y^n_h \Vert_{L^{\infty}}^2 \Bigr)^{-1}, \qquad n \in \mathbb{N}_0.
\end{equation}
\end{Lemma}
{\it Proof.} Inserting ${\bf \chi} = \vu^n_h$ into (\ref{eq:flowd}) we obtain after some calculations
\begin{eqnarray}
\lefteqn{ \frac{1}{\tau_n} \Vert \vu^{n+1}_h - \vu^n_h \Vert^2 + \frac{\sigma \epsilon}{2} \Vert D ( \vu^{n+1}_h - \vu^n_h) \Vert^2 +
\frac{\sigma}{2 \epsilon} \Vert   \vu^{n+1}_h - \vu^n_h \Vert^2 } \nonumber \\
& & + \sigma  \int_{\Omega} \bigl( \frac{\epsilon}{2} | D \vu^{n+1}_h |^2 + \frac{1}{2 \epsilon} 
(1 - | \vu^{n+1}_h |^2) \bigr) dx
- \sigma \int_{\Omega} \bigl( \frac{\epsilon}{2} | D \vu^n_h |^2 + \frac{1}{2 \epsilon} ( 1 - | \vu^n_h |^2) \bigr) dx  
\nonumber \\[2mm]
& \leq & 
\int_{\Omega}  a(\vu^{n+1}_h)  \nabla  y_h^n \cdot   \nabla p_h^n dx - 
\int_{\Omega}   a(\vu^n_h)  \nabla  y_h^n \cdot   \nabla p_h^n dx  \equiv  :  I + II. \label{energ1}  
\end{eqnarray} 
Using (\ref{eq:stated}) for $y^n_h$ and $y^{n+1}_h$ with test function $p^n_h$ as well as (\ref{eq:duald}) we may
rewrite $II$ as follows:
\begin{eqnarray}
II & = &  - \int_{\Omega} a(\vu^{n+1}_h) \nabla y^{n+1}_h \cdot  \nabla p_h^n dx \label{energ2}  \\
& = & - \int_{\Omega} a(\vu^{n+1}_h) \nabla y^{n+1}_h \cdot  \nabla p_h^{n+1} dx +  \int_{\Omega} a(\vu^{n+1}_h) \nabla y^{n+1}_h 
\cdot  \nabla (p^{n+1}_h - p_h^n) dx \nonumber \\
& = & -  ( y^{n+1}_h - y^h_{obs},y^{n+1}_h)_{\mathcal O} +  \int_{\Omega} a(\vu^{n+1}_h) \nabla y^{n+1}_h 
\cdot  \nabla (p^{n+1}_h - p_h^n) dx \equiv II_1 + II_2. \nonumber
\end{eqnarray}
Using again (\ref{eq:duald}) we may write
\begin{eqnarray*}
II_1 & = &  - \frac{1}{2}  || y^{n+1}_h - y^h_{obs} ||_{\mathcal O}^2  +  \frac{1}{2}|| y^n_h - y^h_{obs} ||_{\mathcal O}^2
- \frac{1}{2}  || y^{n+1}_h - y^n_h ||_{\mathcal O}^2   - ( y^{n+1}_h - y^h_{obs},y^n_h)_{\mathcal O}  \\
& = &  - \frac{1}{2}  || y^{n+1}_h - y^h_{obs} ||_{\mathcal O}^2  +  \frac{1}{2}|| y^n_h - y^h_{obs} ||_{\mathcal O}^2
- \frac{1}{2}  || y^{n+1}_h - y^n_h ||_{\mathcal O}^2  \\
& &  - \int_{\Omega} a(\vu^{n+1}_h) \nabla y^n_h \cdot \nabla p^{n+1}_h dx,
\end{eqnarray*}
while
\begin{displaymath}
II_2 = \int_{\Omega} a(\vu^n_h) \nabla y^n_h  \cdot  \nabla (p^{n+1}_h - p_h^n) dx.
\end{displaymath}
Inserting the above identities  into (\ref{energ2}) and combining it with (\ref{energ1}) we obtain
\begin{eqnarray}
\lefteqn{ \hspace{-2cm}
 \bigl( \frac{1}{\tau_n} + \frac{\sigma}{2 \epsilon} \bigr)  \Vert \vu^{n+1}_h - \vu^n_h \Vert^2 
+ \frac{\sigma \epsilon}{2} \Vert D ( \vu^{n+1}_h - \vu^n_h) \Vert^2  + \frac{1}{2}  || y^{n+1}_h - y^n_h ||_{\mathcal O}^2 +  J_{\epsilon,h}(\vu^{n+1}_h) }    \nonumber \\
& \leq &   J_{\epsilon,h}(\vu^n_h) + \int_{\Omega} \bigl( a(\vu^n_h) - a(\vu^{n+1}_h) \bigr) \nabla y^n_h 
\cdot \nabla (p^{n+1}_h - p^n_h) dx \nonumber \\
& \leq & J_{\epsilon,h}(\vu^n_h) + \hat{a} \Vert \nabla y^n_h \Vert_{L^{\infty}} \Vert \vu^{n+1}_h - \vu^n_h \Vert \, 
\Vert \nabla (p^{n+1}_h - p^n_h) \Vert. \label{energ3}
\end{eqnarray}
It remains to estimate $\Vert \nabla (p^{n+1}_h - p^n_h) \Vert$.  To begin, note that
\begin{displaymath}
\int_{\Omega} a(\vu_h^{n+1})\nabla ( p^{n+1}_h - p_h^n) \cdot  \nabla \chi dx  =   \int_{\Omega} (a(\vu_h^n)-a(\vu^{n+1}_h)) 
\nabla p_h^n \cdot  \nabla \chi dx + (y^{n+1}_h - y^n_h, \chi)_{\mathcal O}
\end{displaymath}
for all $\chi \in V_h$. Inserting $\chi = p^{n+1}_h - p_h^n$ we deduce that
 \begin{eqnarray*}
\lefteqn{a_{min}  \Vert \nabla ( p^{n+1}_h - p_h^n) \Vert^2 } \\[2mm]
& \leq & \hat{a} \, \Vert \nabla p^n_h \Vert_{L^{\infty}} \,
  \Vert \vu^{n+1}_h - \vu^n_h \Vert \,  \Vert \nabla (p^{n+1}_h - p^n_h) \Vert + \Vert y^{n+1}_h - y^n_h \Vert_{\mathcal O}
\Vert p^{n+1}_h - p^n_h \Vert_{\mathcal O},
\end{eqnarray*}
 which implies in view of (\ref{poincare1})
\begin{displaymath} 
 \Vert \nabla ( p^{n+1}_h - p_h^n) \Vert \leq \frac{\hat{a}}{a_{min}} \Vert \nabla p^n_h \Vert_{L^{\infty}} \,   
\Vert \vu^{n+1}_h - \vu^n_h \Vert + \frac{1}{\sqrt{\hat{c}}} \frac{1}{a_{min}} \Vert y^{n+1}_h - y^n_h \Vert_{\mathcal O}.
\end{displaymath}
Inserting the above bounds into (\ref{energ3}) and using (\ref{timestepcond}) we infer
\begin{eqnarray*}
\lefteqn{  \frac{1}{\tau_n}  \Vert \vu^{n+1}_h - \vu^n_h \Vert^2 
+ \frac{\sigma \epsilon}{2} \Vert D ( \vu^{n+1}_h - \vu^n_h) \Vert^2  + \frac{1}{2}  || y^{n+1}_h - y^n_h ||_{\mathcal O}^2 +  J_{\epsilon,h}(\vu^{n+1}_h) - J_{\epsilon,h}(\vu^n_h) }     \\
& \leq &  \frac{\hat{a}^2}{a_{min}} \Vert \nabla y^n_h \Vert_{L^{\infty}} \Vert \nabla p^n_h \Vert_{L^{\infty}}
\Vert \vu^{n+1}_h - \vu^n_h \Vert^2 +  \frac{\hat{a}}{a_{min}}\frac{1}{\sqrt{\hat{c}}} \Vert \nabla y^n_h \Vert_{L^{\infty}}
\Vert y^{n+1}_h - y^n_h \Vert_{\mathcal O} \Vert \vu^{n+1}_h - \vu^n_h \Vert \\
& \leq & \Bigl( \frac{\hat{a}^2}{a_{min}} \Vert \nabla y^n_h \Vert_{L^{\infty}} \Vert \nabla p^n_h \Vert_{L^{\infty}}
+ \frac{\hat{a}^2}{a_{min}^2}\frac{1}{2 \hat{c}} \Vert \nabla y^n_h \Vert_{L^{\infty}}^2 \Bigr) \Vert \vu^{n+1}_h - \vu^n_h \Vert^2
+ \frac{1}{2}  || y^{n+1}_h - y^n_h ||_{\mathcal O}^2 \\
& \leq &  \bigl( \frac{1}{\tau_n} -1 \bigr) \Vert \vu^{n+1}_h - \vu^n_h \Vert^2
+ \frac{1}{2}  || y^{n+1}_h - y^n_h ||_{\mathcal O}^2,
\end{eqnarray*}
and the result follows. \qed

\hspace{3mm}

\begin{Corollary} Let $\vu^0_h \in \mathcal K_h$. Then the time steps $\tau_n$ in (\ref{eq:flowd}) can be chosen in such
a way that $\tau_n \geq \gamma>0, n \in \mathbb{N}$, where $\gamma$ depends on the data and possibly on $h$. For this choice
the sequence $(\vu^n_h)_{n \in \mathbb{N}}$  generated by (\ref{eq:flowd}) has  a subsequence $(\vu^{n_k}_h)_{k \in 
\mathbb{N}}$ such that $\vu^{n_k}_h \rightarrow \vu_h$ in $W^{1,\infty}(\Omega,\mathbb{R}^r), k \rightarrow
\infty$ and $\vu_h$ satisfies (\ref{discvar1}).
\end{Corollary}
{\it Proof.} Lemma \ref{monotonicity} implies that 
\begin{displaymath}
 \sum_{n=0}^{\infty} \Vert \vu^{n+1}_h - \vu^n_h \Vert^2  \leq  J_{\epsilon,h}(\vu^0_h), \quad   
\sup_{n \in \mathbb{N}_0} J_{\epsilon,h}(\vu^n_h)  \leq  J_{\epsilon,h}(\vu^0_h),
\end{displaymath}
so that $(\vu^n_h)_{n \in \mathbb{N}}$ is bounded in $H^1(\Omega,\mathbb{R}^r)$ and 
\begin{equation} \label{difzero}
\displaystyle
\lim_{n \rightarrow \infty} \Vert \vu^{n+1}_h - \vu^n_h \Vert =0. 
\end{equation}
In addition we infer from (\ref{eststated}) and
(\ref{eq:duald}) that $(y^n_h)_{n \in \mathbb{N}}$ and $(p^n_h)_{n \in \mathbb{N}}$ are also bounded in $H^1(\Omega)$ and
hence also in $W^{1,\infty}(\Omega)$  since $\mbox{dim}V_h < \infty$. In particular, we infer from (\ref{timestepcond})
that the time steps $\tau_n$ can be chosen to be bounded from below by a positive constant. As a result
there exists a subsequence $(\vu^{n_k}_h,y^{n_k}_h,p^{n_k}_h)_{k \in \mathbb{N}}$
and $(\vu_h,y_h,p_h) \in \mathcal K_h \times V_h \times V_h$ such that
\begin{displaymath}
\vu^{n_k}_h \rightarrow \vu_h \mbox{ in } W^{1,\infty}(\Omega,\mathbb{R}^r), \quad y^{n_k}_h \rightarrow y_h, 
\quad p^{n_k}_h \rightarrow p_h 
\quad \mbox{ in } W^{1,\infty}(\Omega) \mbox{ and a.e. in } \Omega.
\end{displaymath}
In particular, $y_h= S_h(\vu_h)$ and $p_h$ satisfies (\ref{eq:duald}). We finally deduce from (\ref{eq:flowd})
\begin{eqnarray*}
\lefteqn{   \sigma  \int_{\Omega} \bigl( \eps D \vu^{n_k+1}_h \cdot  D (\vv_h- \vu^{n_k+1}_h) -  
 \frac1{\eps} \vu^{n_k}_h \cdot (\vv_h - \vu^{n_k+1}_h) \bigr) dx } \\
& & -  \int_{\Omega} (a(\vv_h)-a(\vu^{n_k+1}_h)) \, \nabla y^{n_k}_h \cdot  \nabla p^{n_k}_h dx \geq
 -\frac{1}{\tau_{n_k}} \int_{\Omega} (\vu^{n_k+1}_h - \vu^{n_k}_h) \cdot (\vv_h - \vu^{n_k+1}_h) dx
\end{eqnarray*} 
for all $\vv_h \in \mathcal K_h$.
Recalling (\ref{difzero}) as well as $\tau_{n_k} \geq \gamma$  we find that $\vu_h$ is a solution of (\ref{discvar1}) 
by passing to the limit $k \rightarrow \infty$. \qed

\setcounter{equation}{0}
\section{Computational examples}\label{s:ce}

We use a preconditioned biconjugate gradient stabilized solver for the stationary forward problem (\ref{eq:stated}) and the adjoint problem (\ref{eq:duald}). 
To solve (\ref{eq:flowd}) we use the primal-dual active set method presented in \cite{BlaGarSar13}, where 
the resulting system of linear equations is solved by 
applying the direct solver UMFPACK \cite{Dav07}.


We set 
\begin{equation}
y_{obs} = \tilde{y}_h+\Lambda n(x),
\label{eq:noise}
\end{equation}
where $n(x)$ is a random variable with the standard normal zero mean distribution, 
$\Lambda\in\mathbb{R}$ and $\tilde{y}_h$ is the solution of 
$$
\int_{\Omega} a(\tilde{u}_h)\nabla \tilde{y}_{h} \cdot  \nabla \chi dx = \int_{\partial \Omega} g_h\chi do 
~~\forall\chi \in V_h
$$
where $\tilde{u}_h$ defines the objective curve.  

There is one regularisation parameter $\sigma$. When the data is noisy we expect that a suitable size of $\sigma$ is obtained by balancing the fidelity term with the regularisation term in the objective functional. The size of $\varepsilon$ is determined by the need to obtain an accurate approximation of the regularised problem. We note that the thickness of the interfacial layer between bulk regions is proportional to $\varepsilon$. 
In order to resolve this interfacial layer we need to choose $h\ll \varepsilon$, see \cite{DecDziEll05} for details. 
Typically reasonable results are obtained with around $8$ to $10$ elements across the interface.  Away from the interface
$h$ can be chosen larger and hence adaptivity in space can heavily
speed up computations. In fact we use the finite element toolbox
Alberta 2.0, see \cite{SchSie05}, for adaptivity
and we implemented the same mesh refinement strategy as in 
\cite{BarNurSty04}, i.e. a fine mesh is constructed for all variables
$\mathbf{u}_h^{n+1}, {y}_h^n $ and ${p}_h^n $ where
$ 0 < (u_h^{n})_i < 1$ for at least one index $i\in\{1,\ldots,r\}$ and with a coarser mesh present in the bulk regions where $(u_h^{n})_i = 0$ 
or $(u_h^{n})_i = 1$ for all $i \in \{1,\ldots, r\}$. In Figure \ref{f:mesh} we display a plot of the triangulation of $\Omega$ which 
illustrates the finer mesh within the interface.

In our computations we found it convenient to choose $h_{min} = \tfrac1{256}$ as the minimal diameter, $h_{max} = \tfrac1{64}$ as  
the maximal diameter of all elements and we set {\color{black}$\tau_n = 0.01/\varepsilon$}. 
{\color{black} The stopping criteria we used to terminate the algorithm was the size of the residual to the 
first order optimality condition, i.e. $\|(\vu_h^{n+1}-\vu_h^n)/\tau_n\|\leq 1.0e^{-3}$. 
For each computation we state the number of iterations, $L$, required to reach this stopping criteria.} 

In the case $r=2$ we have  $u_2=1 -u_1$ and the vector-valued Allen-Cahn inequality with two order parameters is reduced 
in the computations to a scalar Allen-Cahn inequality. 
\begin{figure}[htbp]
\begin{center}
\includegraphics[width=.44\textwidth,angle=0]{{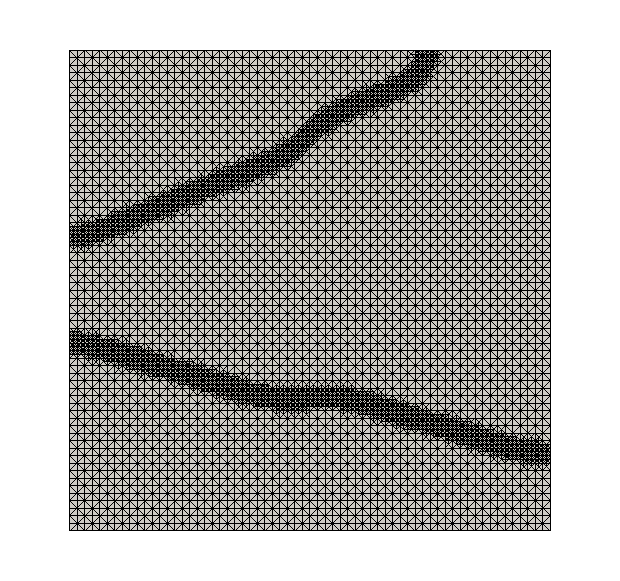}}\hspace{0mm}

\caption{A converged triangulation}
\label{f:mesh}
\end{center}
\end{figure}


\subsection{Results with $r=2$ and $d=2$}
In this section we see how our method compares with the one presented in \cite{ItoKunLi01}.  
In all the computations unless otherwise stated we set $\Omega=(-1,1)^2$, 
$J_{fid}(\Gamma):= ||y_{\Gamma}-y_{obs}||_{L^2(\Omega)}^{2}$, 
$\displaystyle{\varepsilon=\frac1{16\pi}}$, $a_1=3$, $a_2=0.5$, $\sigma=0.0001$, $\Lambda=0.05$ and 
$$
g_h(x,y)=\left\{ \begin{array}{cl} 
-0.5&\mbox{ if }x=-1\mbox{ or } y=-1\\
0.5&\mbox{ if }x=1\mbox{ or } y=1.
\end{array}\right.
$$

\vspace{1mm}

Figure \ref{f1} displays the results we obtain when using the same initial curve (a circle of radius $0.6$) and objective curve (a `skinny' ellitpse, 
$x^2/(0.07)^2 + y^2/(0.5)^2=1$) that are used in Section 4.1 of  \cite{ItoKunLi01}. In this simulation we set $\Lambda=0$, as in  \cite{ItoKunLi01}. 
The left hand plot in Figure \ref{f1} displays the initial curve, the centre plot the objective curve and the right hand plot the 
computed solution $\mathbf{u}_h^n$. {\color{black} The number of iterations required to reach the stopping criteria was $L=4417$. }

Figure \ref{f2} takes the same form as Figure \ref{f1} except that this time we compare our results with those displayed in Section 4.4 of  \cite{ItoKunLi01}. 
The initial curve is again a circle of radius $0.6$ while the objective curve consists of two objects
$$
\frac{(x+0.35)^2}{(0.25)^2} + \frac{(y+0.35)^2}{(0.3)^2}=1~~\mbox{and}~~ \frac{(x-0.35)^2}{(0.2)^2} + \frac{(y-0.35)^2}{(0.2)^2}=1,
$$
as in \cite{ItoKunLi01} we set $\Lambda=0$.  {\color{black} The number of iterations required to reach the stopping criteria was $L=11117$. }
From this example we see that our phase field model successfully deals with topological change. 

{\color{black}In Figure \ref{fe1} we plot the Residual $:= \|(\vu_h^{n+1}-\vu_h^n)/\tau_n\|$, 
$J^{fid}_{\varepsilon,h}(\mathbf{u}_h):=\frac12||S_{h}(\vu_{h})-y^{h}_{obs}||^2_{\mathcal O}$, 
$J^{reg}_{\varepsilon,h}(\mathbf{u}_h):=\sigma  \int_\Omega \bigl( \frac{\eps}2| D \vu_h|^2+\frac1{2\eps} (1 - | \vu_h |^2) \bigr)dx$ 
and  $J_{\varepsilon,h}(\mathbf{u}_h)$
versus iteration number, for the first $2000$ iterations,}
for the computations displayed in Figures \ref{f1} and \ref{f2}. 
{\color{black} From this figure we see that, in both computations, for the first $50$ iterations there is a steep decrease in 
$J_{\epsilon,h}(\vu_h)$ and after that the decrease is much more gradual. We also see that the Residual decreases at a much slower rate than 
$J_{\epsilon,h}(\vu_h)$.}
{\color{black} In Figure \ref{f:ekk} 
we display two intermediate results from the set-up in Figure \ref{f1}; 
the plots display $\mathbf{u}_h^n$ 
after $150$ iterations (left hand plot), 
after $500$ iterations (centre plot) and after $4417$ iterations, once the iteration has converged (right hand plot). 
From this figure we see that after $500$ iterations the solution is approximating the shape of the objective curve 
reasonably well although the curve is not yet defined by a well defined interfacial region.   }

\begin{figure}[htbp]
\begin{center}
\includegraphics[width=.32\textwidth,angle=0]{{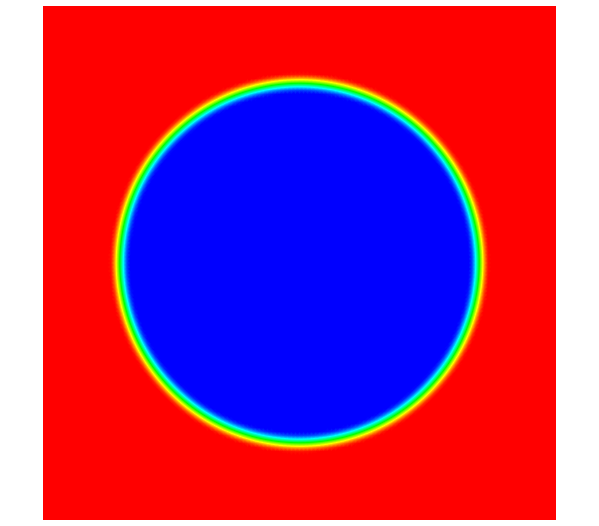}}\hspace{0mm}
\includegraphics[width=.32\textwidth,angle=0]{{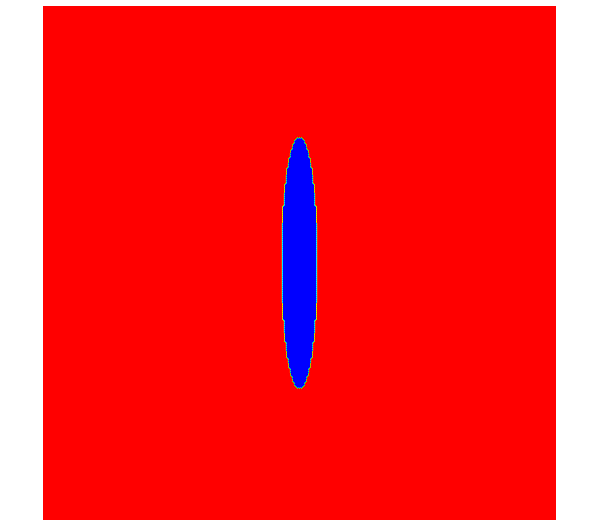}}\hspace{0mm}
\includegraphics[width=.32\textwidth,angle=0]{{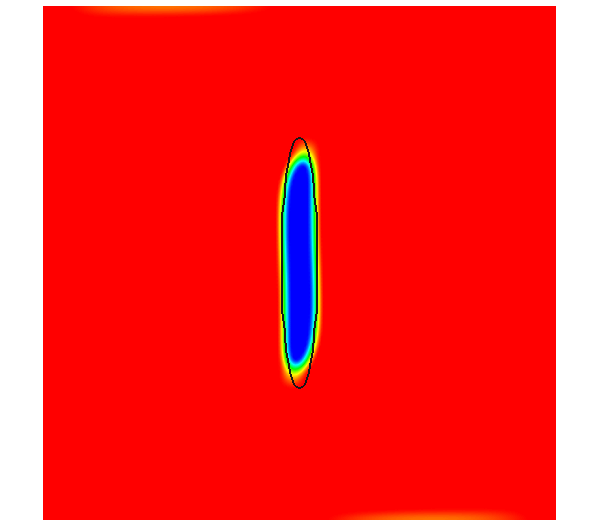}}\hspace{0mm}
\caption{$J_{fid}(\Gamma):= ||y_{\Gamma}-y_{obs}||_{L^2(\Omega)}^{2}$, initial curve (left hand plot), objective curve (centre plot), 
$\mathbf{u}_h^n$ (right hand plot)}
\label{f1}
\end{center}
\end{figure} 

 \begin{figure}[htbp]
\begin{center}
\includegraphics[width=.32\textwidth,angle=0]{{ikl_id_0p7.png}}\hspace{0mm}
\includegraphics[width=.32\textwidth,angle=0]{{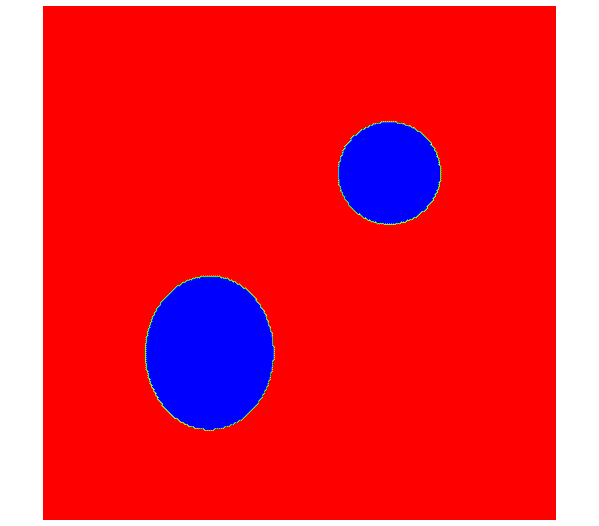}}\hspace{0mm}
\includegraphics[width=.32\textwidth,angle=0]{{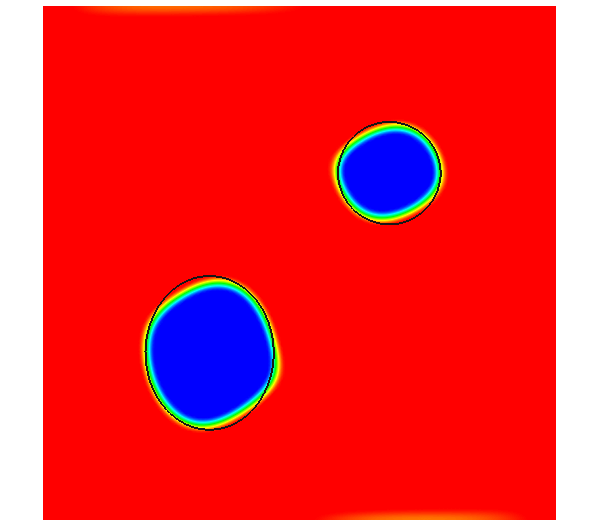}}\hspace{0mm}
\caption{$J_{fid}(\Gamma):= ||y_{\Gamma}-y_{obs}||_{L^2(\Omega)}^{2}$, initial curve (left hand plot), objective curve (centre plot), 
$\mathbf{u}_h^n$ (right hand plot)}
\label{f2}
\end{center}
\end{figure}


 \begin{figure}[htbp]
\begin{center}
\includegraphics[width=.34\textwidth,angle=0]{{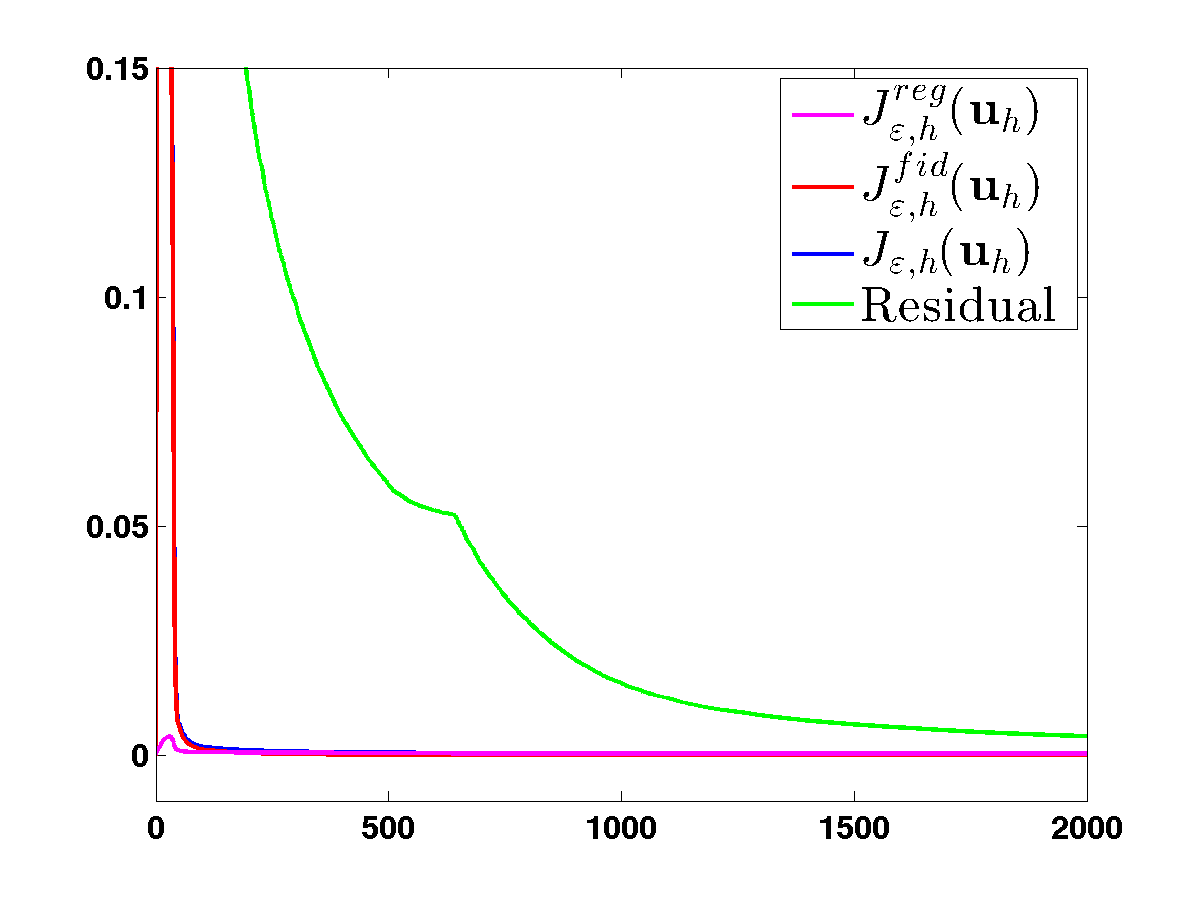}}\hspace{-4mm}
\includegraphics[width=.34\textwidth,angle=0]{{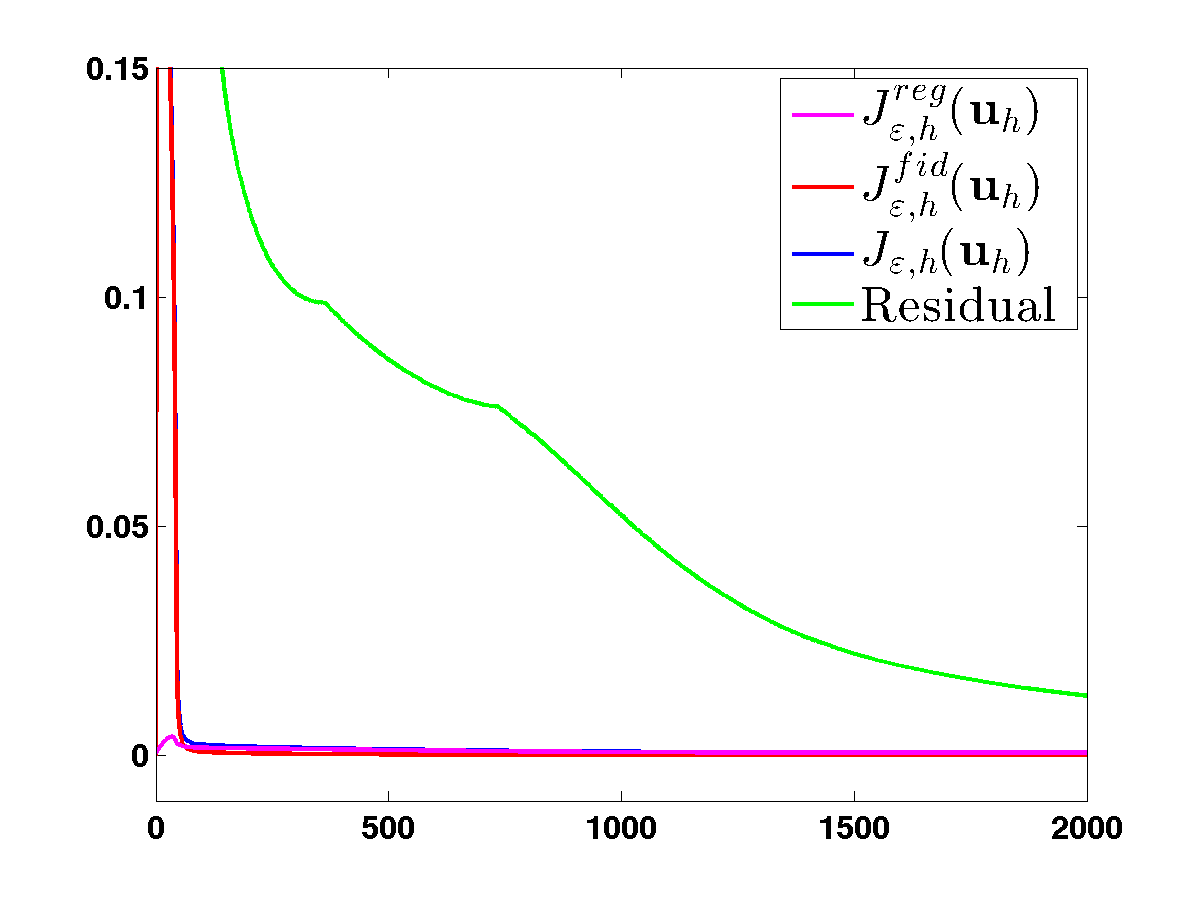}}\hspace{-4mm}
\caption{Plot of $J^{fid}_{\varepsilon,h}(\mathbf{u}_h)$, $J^{reg}_{\varepsilon,h}(\mathbf{u}_h)$, $J_{\varepsilon,h}(\mathbf{u}_h)$ and the Residual, versus the number of iterations: results in Figure \ref{f1} (left plot), results in Figure \ref{f2} (right plot)}
\label{fe1}
\end{center}
\end{figure} 

\begin{figure}[htbp]
\begin{center}
\includegraphics[width=.3\textwidth,angle=0]{{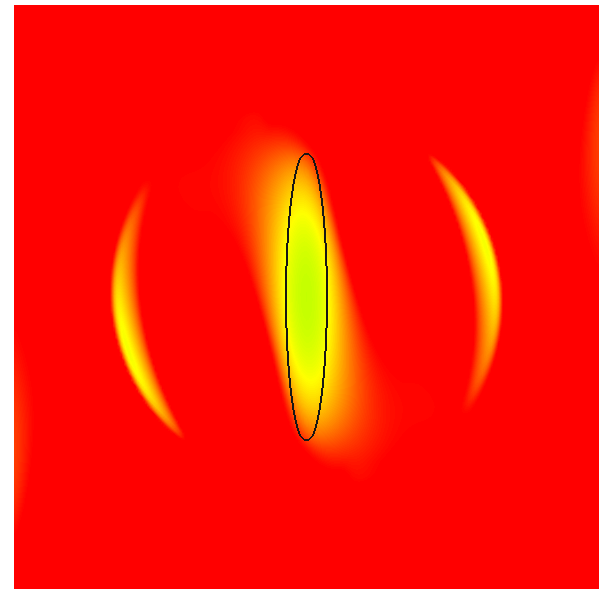}}\hspace{0mm}
\includegraphics[width=.3\textwidth,angle=0]{{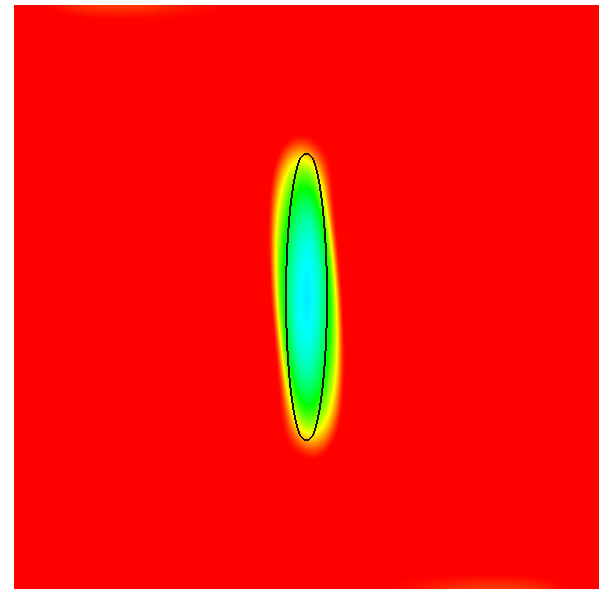}}\hspace{0mm}
\includegraphics[width=.3\textwidth,angle=0]{{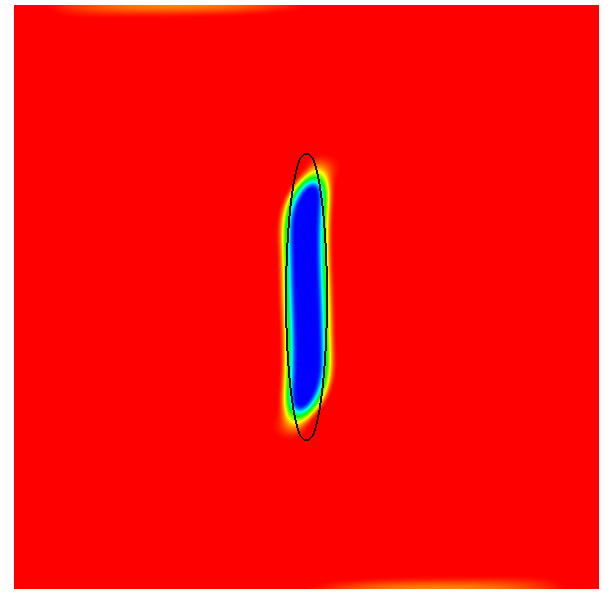}}\hspace{0mm}
\caption{$J_{fid}(\Gamma):= ||y_{\Gamma}-y_{obs}||_{L^2(\Omega)}^{2}$, $\mathbf{u}_h^n$ after $150$ iterations (left plot), $\mathbf{u}_h^n$ after $500$ iterations (centre plot) 
$\mathbf{u}_h^n$ after $4417$ iterations (right plot)}
\label{f:ekk}
\end{center}
\end{figure}

In Figure \ref{f3} we follow the authors in Section 4.2 of \cite{ItoKunLi01} in seeing how noise effects the solution. 
We take the same initial and objective curves as in Figure \ref{f1} and display the solutions obtained with $\Lambda=0.05$ (left hand plot), 
$\Lambda=0.1$ (centre plot)  and $\Lambda=0.2$ (right hand plot).  
{\color{black} The number of iterations required to reach the stopping criteria were $L=4236$, $L=4075$ and $L=8941$ respectively. }

In Figure \ref{f4} we follow the authors in Section 4.5 of \cite{ItoKunLi01} in seeing how the value of the regularisation parameter 
$\sigma$ effects the solution. 
For the initial curve we take a circle of radius $0.7$ and for the objective curve we take the ellipse $x^2/(0.5)^2 + y^2/(0.4)^2=1$. 
For the choice $\Lambda = 0.05$ we display the solutions obtained with $\sigma=0.01$ (top centre) $\sigma=0.001$ (top right) and $\sigma=0.0001$ (bottom left) 
$\sigma=0.000025$ (bottom centre) $\sigma=0.0000025$ (bottom right). 
{\color{black} The number of iterations required to reach the stopping criteria were $L=3918$, $L=9183$, $L=5550$, $L=8441$ and $L=21228$ respectively. }
From this figure we see that $\sigma =  0.001$ and $\sigma=0.0001$ give the best approximations to the objective curve. 
%
 

In Figure \ref{fe2} we plot $J^{fid}_{\epsilon,h}(\vu_h)$, $J^{reg}_{\epsilon,h}(\vu_h)$, $J_{\epsilon,h}(\vu_h)$ and the Residual 
for the first $4000$ iterations, for the 
computations displayed in Figure \ref{f4} with $\sigma =  0.001$, $\sigma=0.0001$ and $\sigma=0.000025$. 
From this figure we see that for $\sigma =  0.001$ the initial decrease in $J_{\epsilon,h}(\vu_h)$ is more gradual than for $\sigma=0.0001$ and $\sigma=0.000025$. 


 \begin{figure}[htbp]
\begin{center}
\includegraphics[width=.32\textwidth,angle=0]{{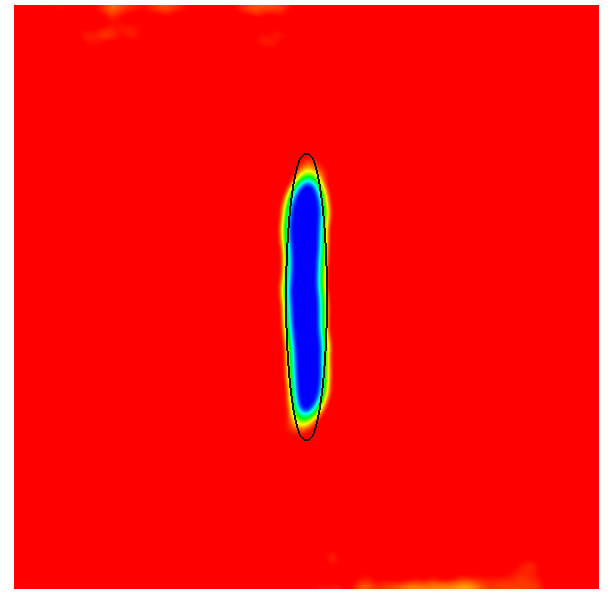}}\hspace{0mm}
\includegraphics[width=.32\textwidth,angle=0]{{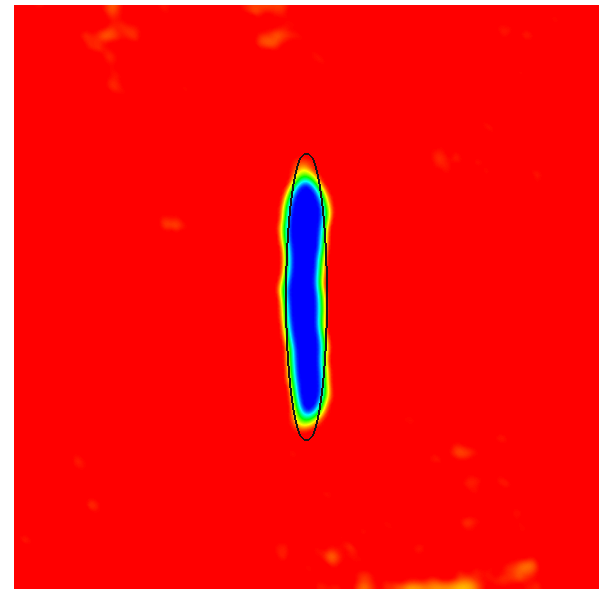}}\hspace{0mm}
\includegraphics[width=.32\textwidth,angle=0]{{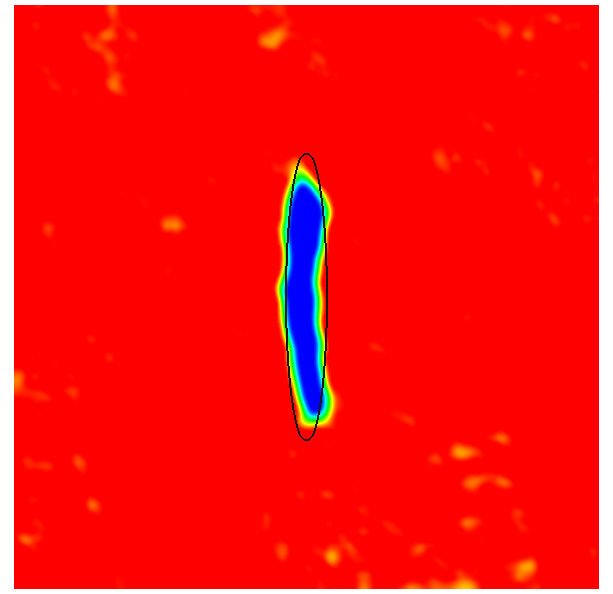}}\hspace{0mm}
\caption{$J_{fid}(\Gamma):= ||y_{\Gamma}-y_{obs}||_{L^2(\Omega)}^{2}$, $\mathbf{u}_h^n$ obtained by taking 
$\Lambda=0.05$ (left hand plot) $\Lambda=0.1$ (centre plot) and $\Lambda=0.2$ (right hand plot)}
\label{f3}
\end{center}
\end{figure} 


\begin{figure}[htbp]
\begin{center}
\includegraphics[width=.26\textwidth,angle=0]{{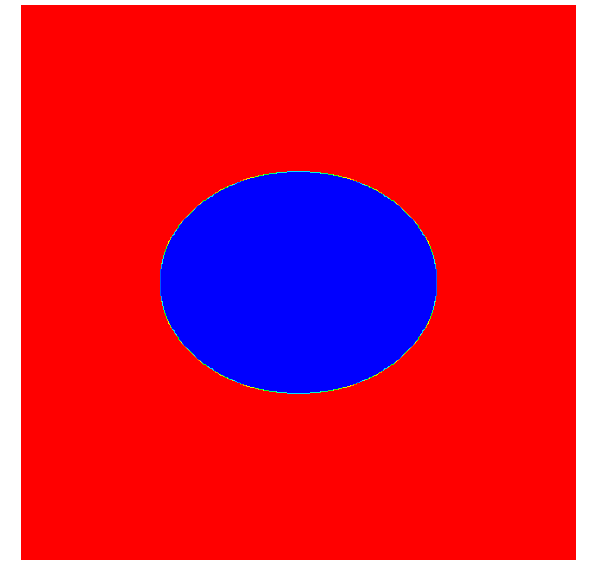}}\hspace{0mm}
\includegraphics[width=.26\textwidth,angle=0]{{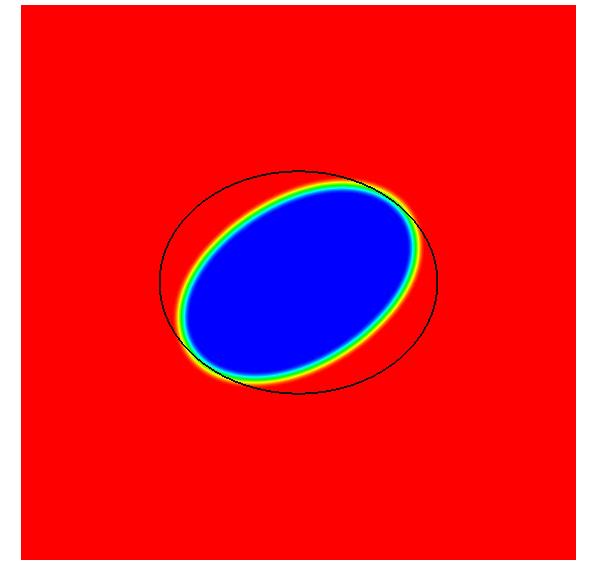}}\hspace{0mm}
\includegraphics[width=.26\textwidth,angle=0]{{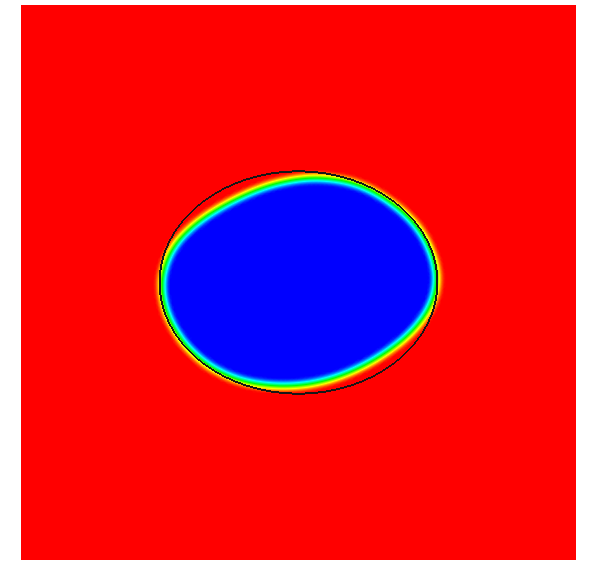}}\hspace{0mm}\\
\includegraphics[width=.26\textwidth,angle=0]{{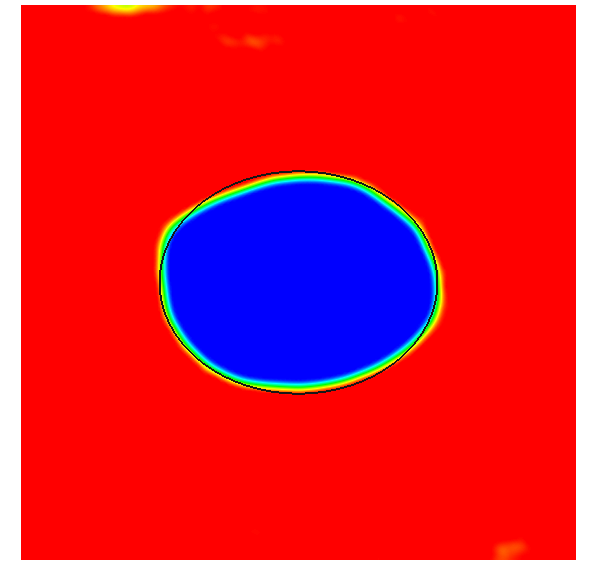}}\hspace{0mm}
\includegraphics[width=.26\textwidth,angle=0]{{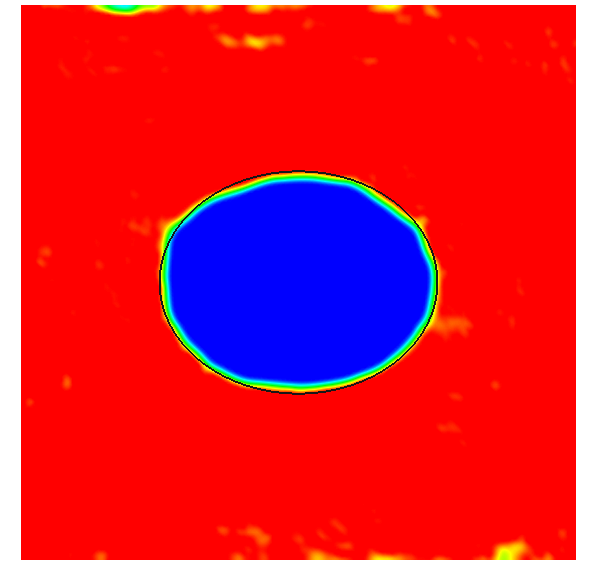}}\hspace{0mm}
\includegraphics[width=.26\textwidth,angle=0]{{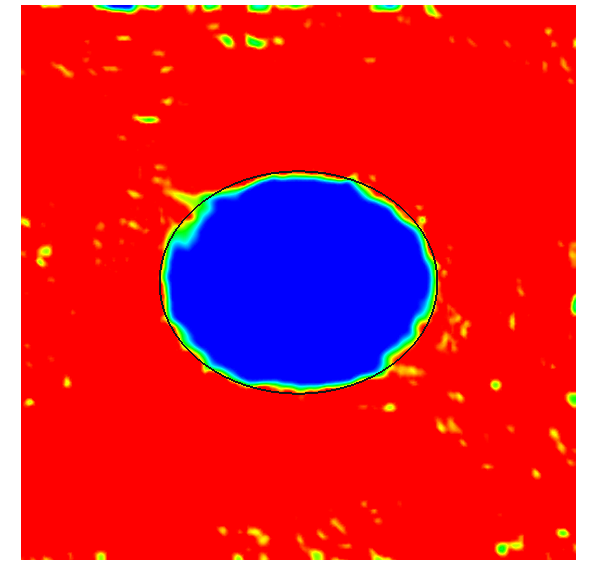}}\hspace{0mm}
\caption{$J_{fid}(\Gamma):= ||y_{\Gamma}-y_{obs}||_{L^2(\Omega)}^{2}$, objective curve (top left), $\mathbf{u}_h^n$ obtained by taking 
$\sigma=0.01$ (top centre) $\sigma=0.001$ (top right) and $\sigma=0.0001$ (bottom left) $\sigma=0.000025$ (bottom centre) $\sigma=0.0000025$ (bottom right)}
\label{f4}
\end{center}
\end{figure} 


 \begin{figure}[htbp]
\begin{center}
\includegraphics[width=.34\textwidth,angle=0]{{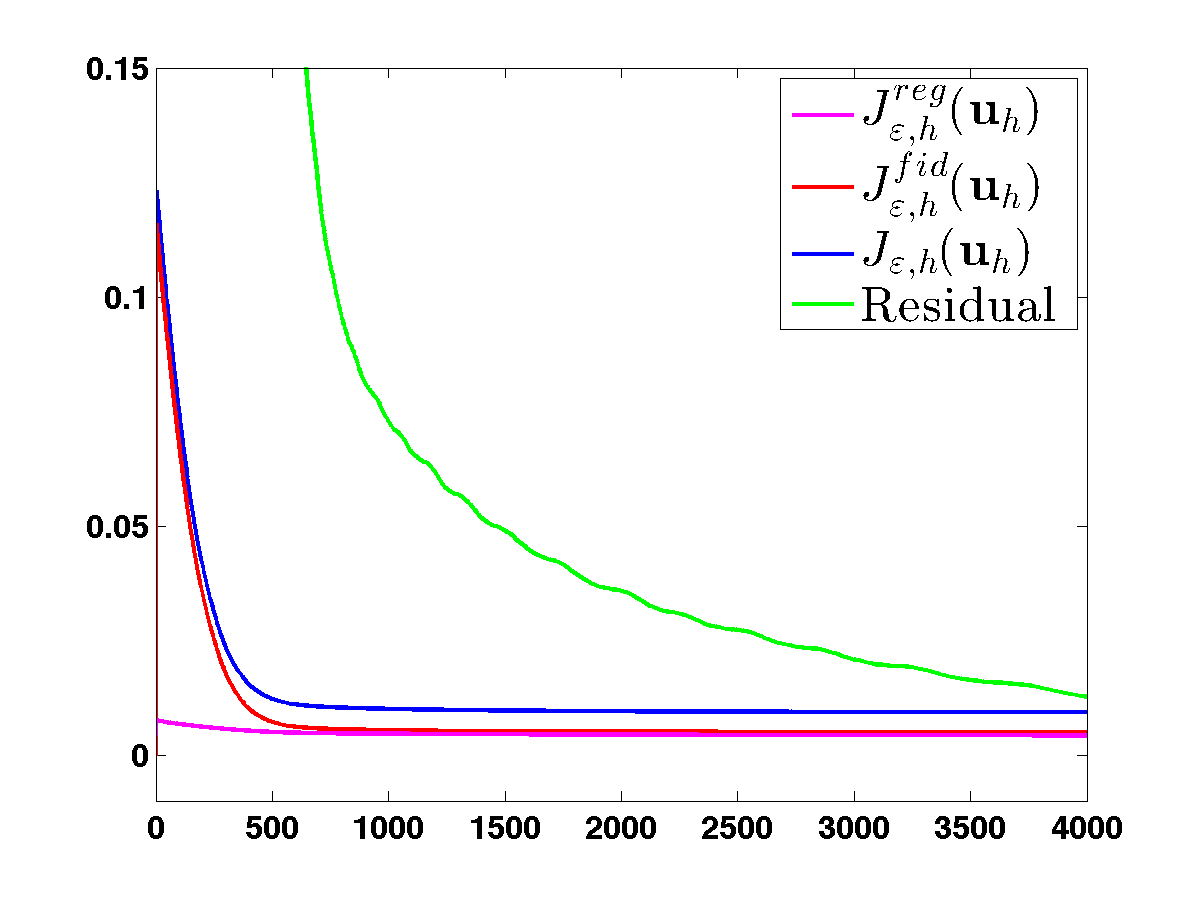}}\hspace{-4mm}
\includegraphics[width=.34\textwidth,angle=0]{{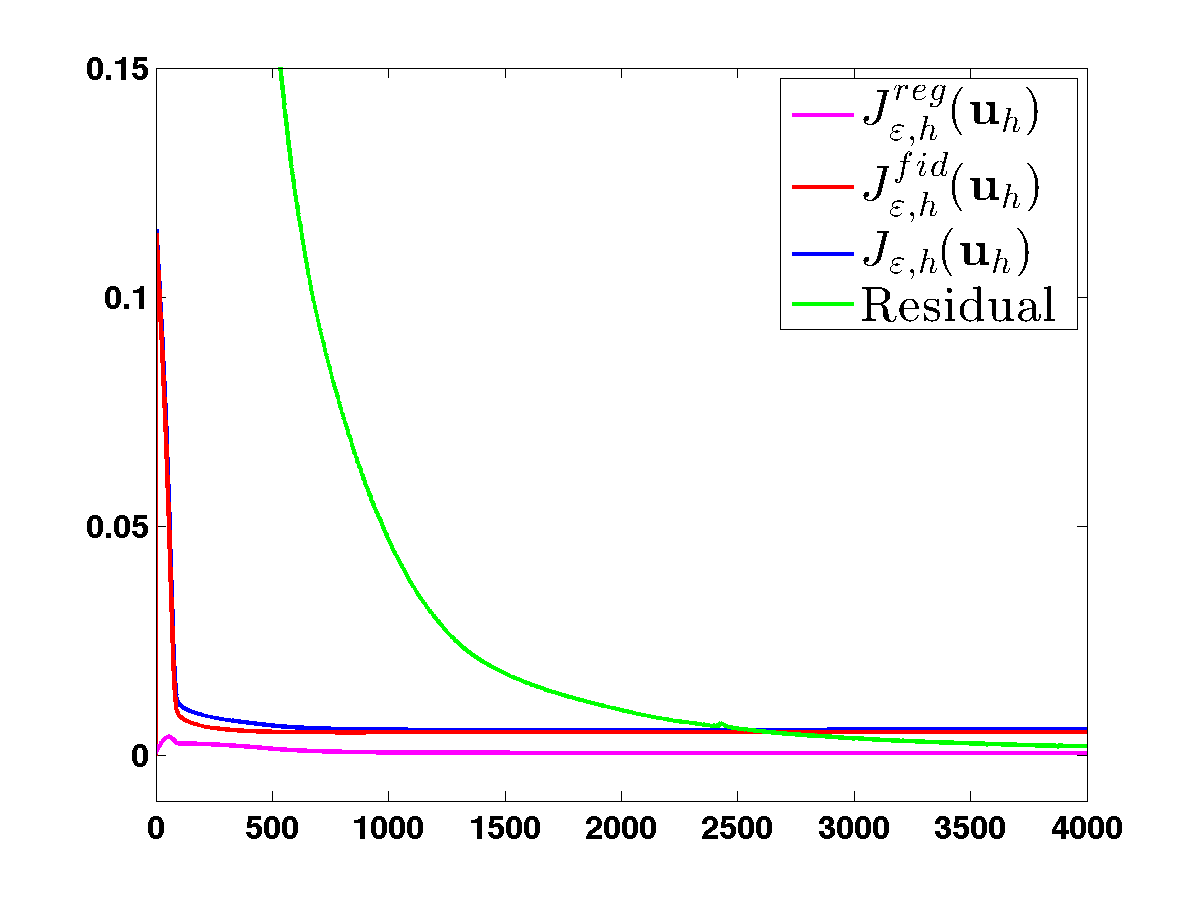}}\hspace{-4mm}
\includegraphics[width=.34\textwidth,angle=0]{{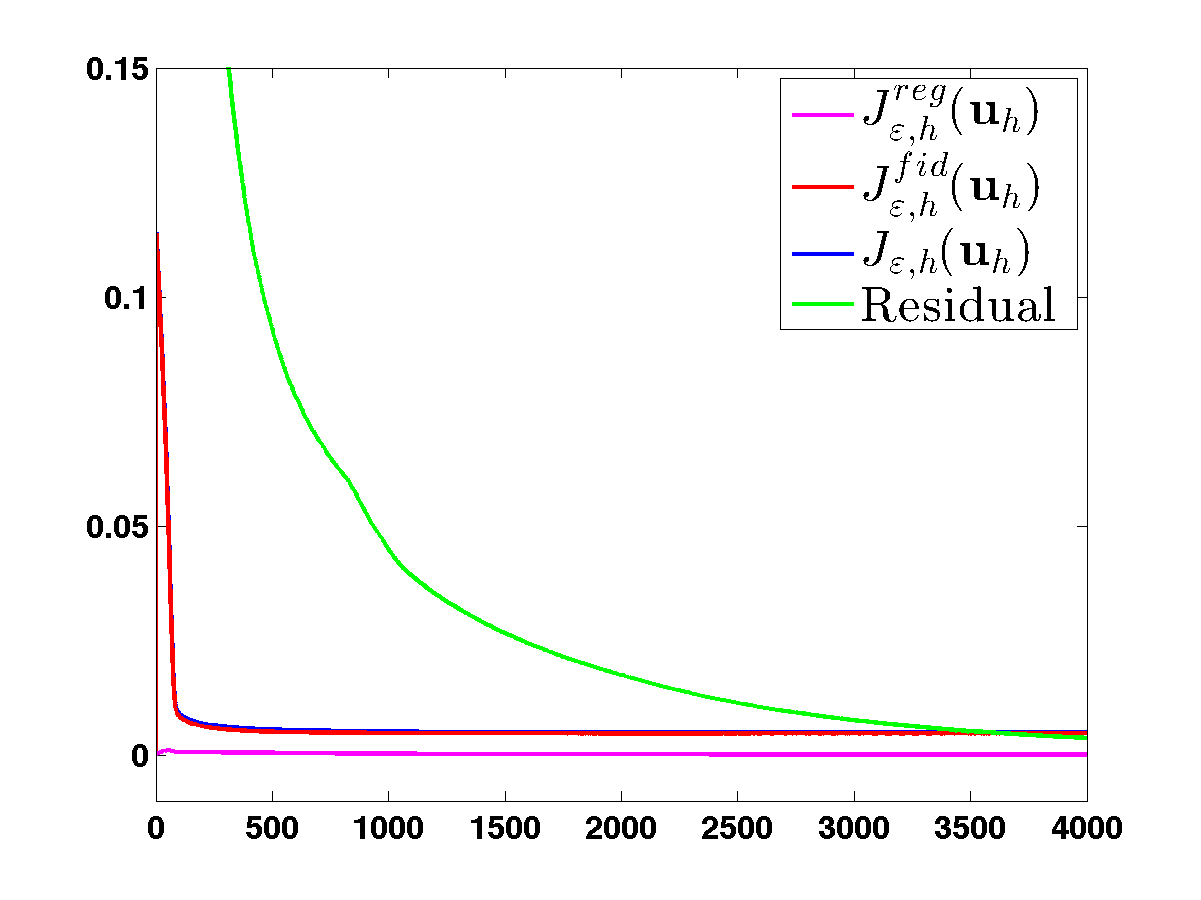}}\hspace{0mm}
\caption{Plot of $J^{fid}_{\varepsilon,h}(\mathbf{u}_h)$, $J^{reg}_{\varepsilon,h}(\mathbf{u}_h)$, $J_{\varepsilon,h}(\mathbf{u}_h)$ and the Residual,
versus the number of iterations for the results in Figure \ref{f4}, $\sigma=0.001$ (left plot) 
$\sigma=0.0001$ (centre plot), $\sigma=0.000025$ (right plot)}
\label{fe2}
\end{center}
\end{figure}

\vspace{2mm}

In Figure \ref{f:0a} we show the effect that the size of $|a_1-a_2|$ has on the solution $\mathbf{u}_h^n$. 
We display the objective curve in the left hand plot and in the subsequent plots we display a zoomed in 
image of the approximate solution, $\mathbf{u}_h^n$, 
at the end of the simulation obtained from decreasing value of $|a_1-a_2|$. We take $a_2=0.5$ in all plots 
and $a_1=1,~3,~7$ in the second, third and fourth plots respectively. 
We see that the approximation to the objective curve improves when $|a_1-a_2|$ increases. 
{\color{black} The number of iterations required to reach the stopping criteria were $L=11891$, $L=5550$ and $L=17072$ respectively. }

\begin{figure}[htbp]
\begin{center}
\includegraphics[width=.24\textwidth,angle=0]{{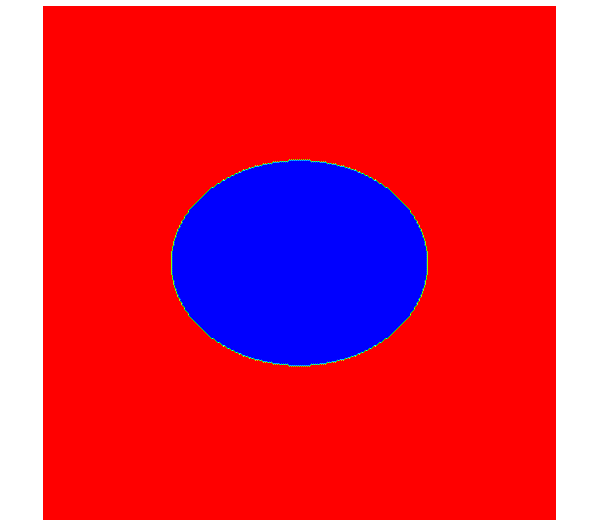}}\hspace{0mm}
\includegraphics[width=.24\textwidth,angle=0]{{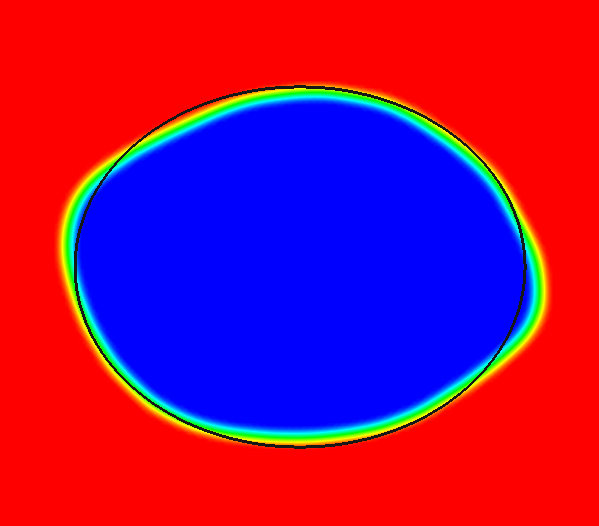}}\hspace{0mm}
\includegraphics[width=.24\textwidth,angle=0]{{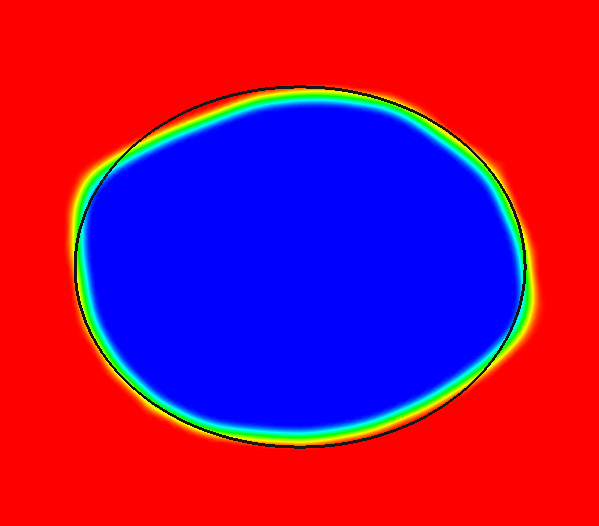}}\hspace{0mm}
\includegraphics[width=.24\textwidth,angle=0]{{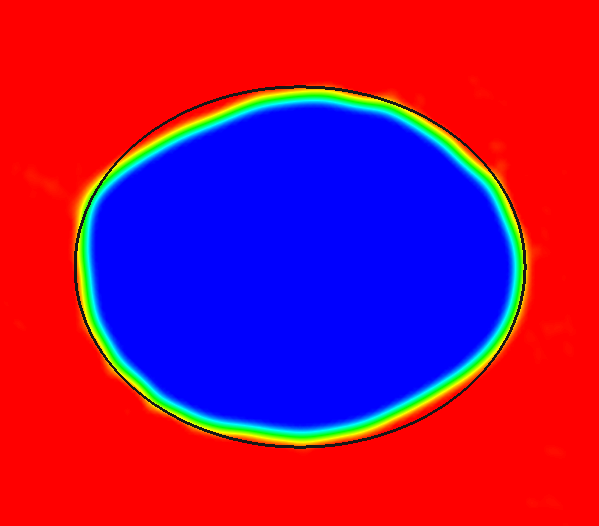}}\hspace{0mm}
\caption{$J_{fid}(\Gamma):= ||y_{\Gamma}-y_{obs}||_{L^2(\Omega)}^{2}$, objective curve (first plot), zoomed in 
plot of $\mathbf{u}_h^n$ obtained by taking $(a_1,a_2)=(1,0.5)$  (second plot), $(a_1,a_2)=(3,0.5)$  
(third plot) and $(a_1,a_2)=(7,0.5)$  (fourth plot)}
\label{f:0a}
\end{center}
\end{figure} 

In Figure \ref{f:1noise} we show the effect that the choice of $\mathcal{O}$ has on the solution $\mathbf{u}_h^n$. 
We compare results obtained by taking 
$J_{fid}(\Gamma):= ||y_{\Gamma}-y_{obs}||_{L^2(\Omega)}^{2}$ to results obtained by taking 
$J_{fid}(\Gamma):= ||y_{\Gamma}-y_{obs}||_{L^2(\partial \Omega)}^{2}$. In these simulations we set $\Lambda = 0.02$. 
We display the objective curve in the left hand plot and the approximate solution $\mathbf{u}_h^n$ 
at the end of the simulation obtained from $J_{fid}(\Gamma):= ||y_{\Gamma}-y_{obs}||_{L^2(\Omega)}^{2}$ (centre plot) 
and $J_{fid}(\Gamma):= ||y_{\Gamma}-y_{obs}||_{L^2(\partial \Omega)}^{2}$ (right plot). From this figure we see that the 
approximation to the objective curve
obtained using $J_{fid}(\Gamma):= ||y_{\Gamma}-y_{obs}||_{L^2(\partial \Omega)}^{2}$ is effected more by the noise than the 
approximation that is obtained using $J_{fid}(\Gamma):= ||y_{\Gamma}-y_{obs}||_{L^2(\Omega)}^{2}$. Furthermore using 
 $J_{fid}(\Gamma):= ||y_{\Gamma}-y_{obs}||_{L^2(\Omega)}^{2}$ gives 
a better approximation to the objective curve than using $J_{fid}(\Gamma):= ||y_{\Gamma}-y_{obs}||_{L^2(\partial \Omega)}^{2}$. 
{\color{black} The number of iterations required to reach the stopping criteria were $L=16151$ and $L=18081$ respectively. }

\begin{figure}[htbp]
\begin{center}
\includegraphics[width=.3\textwidth,angle=0]{{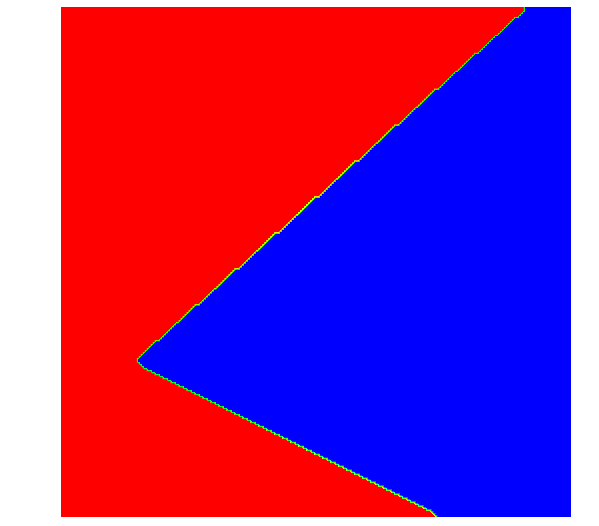}}\hspace{0mm}
\includegraphics[width=.3\textwidth,angle=0]{{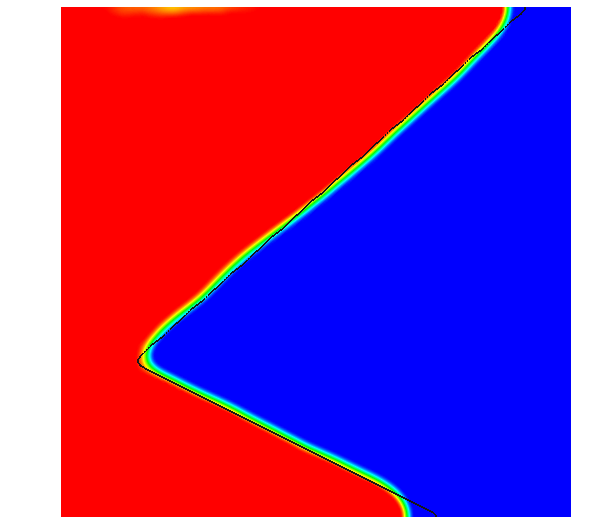}}\hspace{0mm}
\includegraphics[width=.3\textwidth,angle=0]{{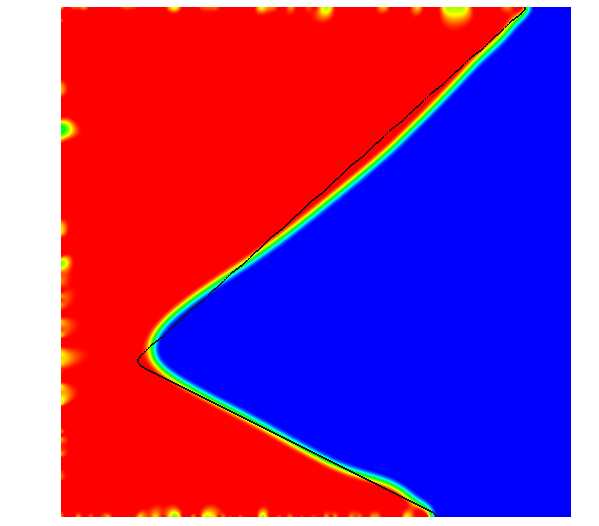}}\hspace{0mm}
\caption{Objective curve (left plot), $\mathbf{u}_h^n$ obtained from $J_{fid}(\Gamma):= ||y_{\Gamma}-y_{obs}||_{L^2(\Omega)}^{2}$ (centre plot) 
and $J_{fid}(\Gamma):= ||y_{\Gamma}-y_{obs}||_{L^2(\partial \Omega)}^{2}$ (right plot) }
\label{f:1noise}
\end{center}
\end{figure}

\vspace{2mm}
%
%
%
%

In Figure \ref{f:1} we display results for three objective curves; 
we plot the objective curves in the upper row and the solution $\mathbf{u}_h^n$ 
at the end of the simulation in the lower row. In these simulations we took $\sigma = 0.00001$ and $\Lambda = 0.005$. 
{\color{black} The number of iterations required to reach the stopping criteria were $L=12720$, $L=22296$ and $L=36036$ respectively. }
From this figure we see that our method results in good approximations of the objective curves.  

\begin{figure}[htbp]
\begin{center}
\includegraphics[width=.24\textwidth,angle=0]{{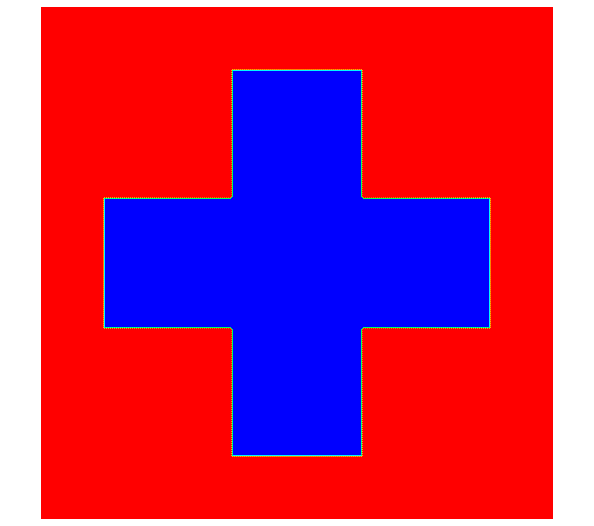}}\hspace{0mm}
\includegraphics[width=.24\textwidth,angle=0]{{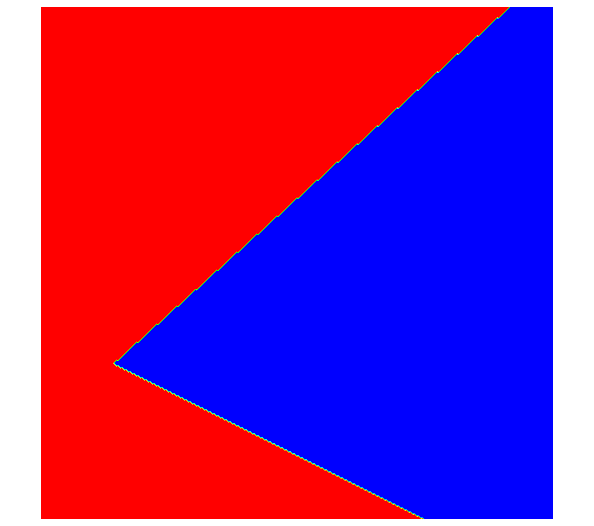}}\hspace{0mm}
\includegraphics[width=.24\textwidth,angle=0]{{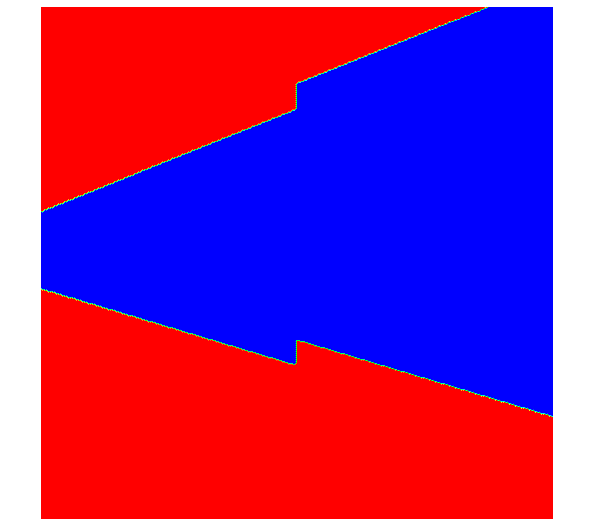}}\\[5mm]
\includegraphics[width=.24\textwidth,angle=0]{{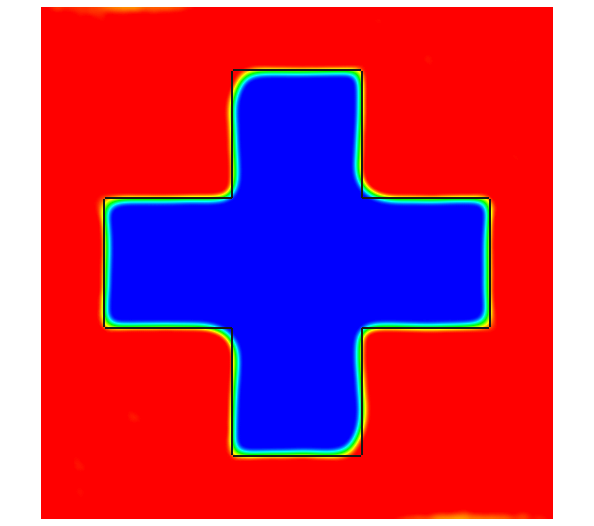}}\hspace{0mm}
\includegraphics[width=.24\textwidth,angle=0]{{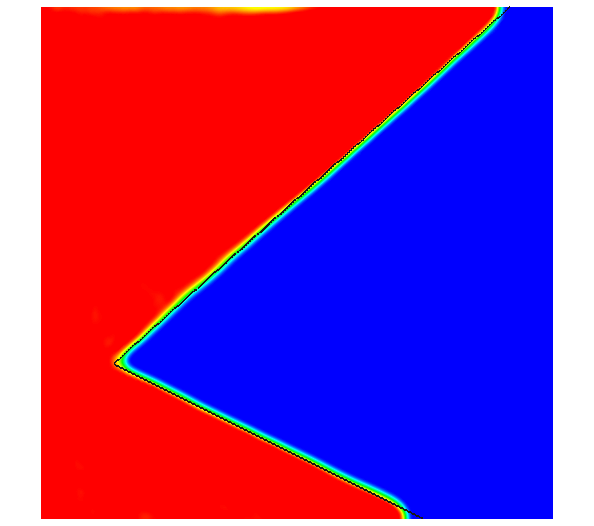}}\hspace{0mm}
\includegraphics[width=.24\textwidth,angle=0]{{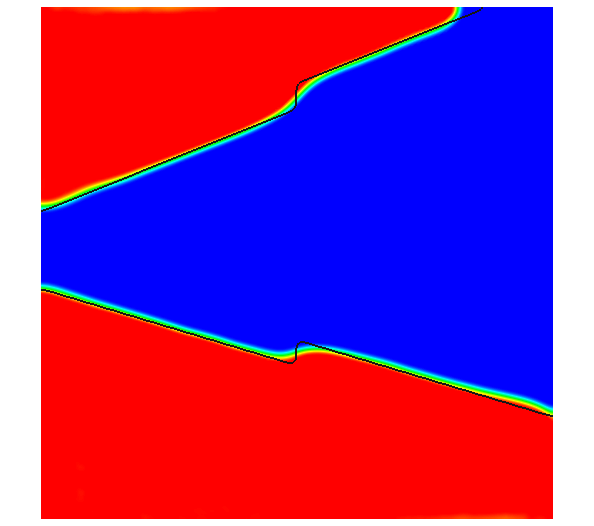}}
\caption{$J_{fid}(\Gamma):= ||y_{\Gamma}-y_{obs}||_{L^2(\Omega)}^{2}$, objective curves (upper plots), $\mathbf{u}_h^n$ (lower plots)}
\label{f:1}
\end{center}
\end{figure}

\subsection{Results with $r=3$ and $d=2$}
In all the computations in this section 
we set $\Omega=(-1,1)^2$, 
$J_{fid}(\Gamma):= ||y_{\Gamma}-y_{obs}||_{L^2(\Omega)}^{2}$, 
$\displaystyle{\varepsilon=\frac1{8\pi}}$, $a_1=0.8$, $a_2=0.2$, $a_3=0.3$, $\sigma=0.001$,  $\Lambda=0.0$ and 
$$
g_h(x,y)=\left\{ \begin{array}{cl} 
0&\mbox{ if }x=\pm 1\\
-0.5&\mbox{ if }y=-1\\
0.5&\mbox{ if }y=1.
\end{array}\right.
$$
In Figure \ref{f:vv1} we display results for four objective curves, for each curve we took random t
initial data for $\mathbf{u}_h^0$.  We plot the objective curves in the upper row and 
the solution $\mathbf{u}_h^n$ at the end of the simulation 
in the lower row. {\color{black} The number of iterations required to reach the stopping criteria were $L=10844$, $L=33574$, $L=31113$ and $L=57373$ respectively. }

\begin{figure}[htbp]
\begin{center}
\includegraphics[width=.24\textwidth,angle=0]{{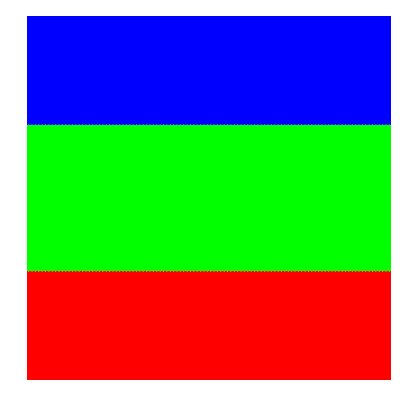}}\hspace{0mm}
\includegraphics[width=.24\textwidth,angle=0]{{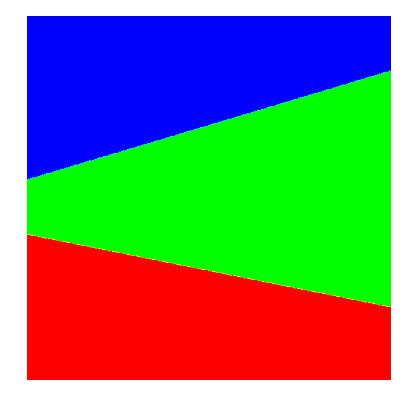}}
\includegraphics[width=.24\textwidth,angle=0]{{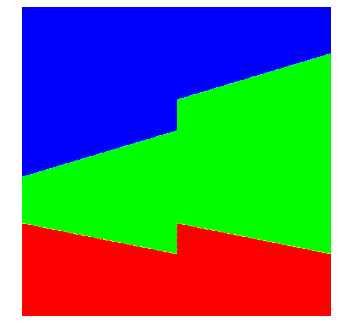}}
\includegraphics[width=.217\textwidth,angle=0]{{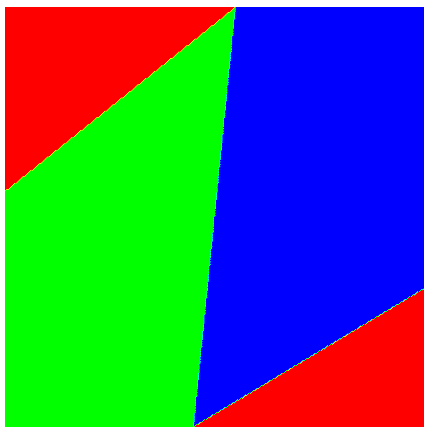}}\\
\includegraphics[width=.24\textwidth,angle=0]{{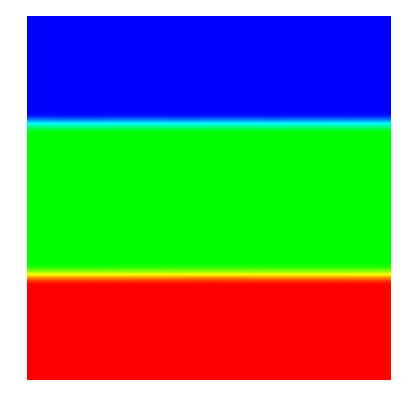}}\hspace{0mm}
\includegraphics[width=.24\textwidth,angle=0]{{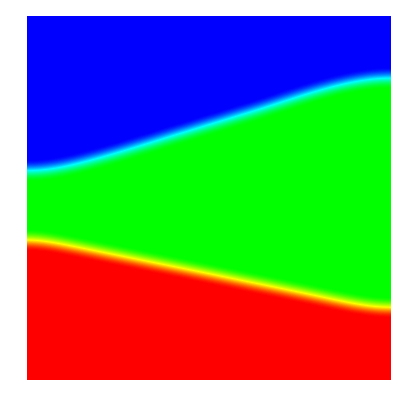}}
\includegraphics[width=.24\textwidth,angle=0]{{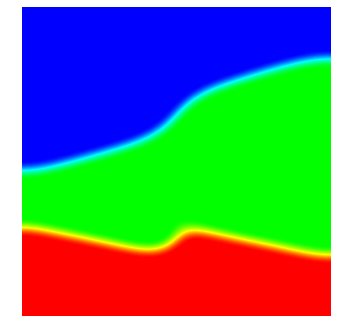}}
\includegraphics[width=.218\textwidth,angle=0]{{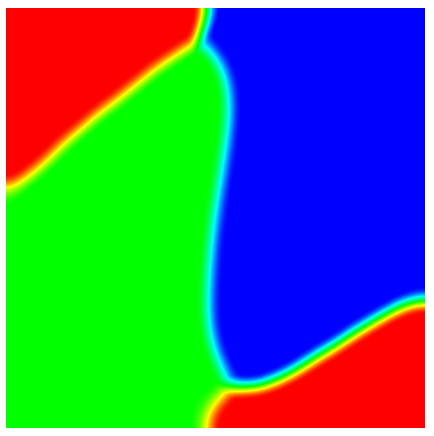}}
\caption{$J_{fid}(\Gamma):= ||y_{\Gamma}-y_{obs}||_{L^2(\Omega)}^{2}$, objective curves (upper plots), $\mathbf{u}_h^n$ (lower plots)}
\label{f:vv1}
\end{center}
\end{figure} 

{\color{black}
\subsection{Summary of the computational results} 
The set-up of the computational examples presented in Figures \ref{f1}, \ref{f2}, \ref{f3} and \ref{f4} are taken from examples presented in \cite{ItoKunLi01}. 
The closeness of the approximated curve to the objective curve in the results that we present in Figures \ref{f3} and \ref{f4} is of a similar order 
to the results presented in \cite{ItoKunLi01}. In the case of Figure \ref{f1} the level set method used in  \cite{ItoKunLi01} yields a 
better approximation to the skinny ellipse than our phase field model while in the case of Figure \ref{f2} our results are a substantial improvement 
on the ones in \cite{ItoKunLi01} as the level set method is unable to deal with the topological change required in this example 
whereas the phase field model successfully deals with it. 
}

\setcounter{equation}{0}
\section{Appendix}

\begin{Theorem} \label{feps}
Let $F_{\epsilon}:X \rightarrow \mathbb{R} \cup \lbrace \infty \rbrace$ be defined by
\begin{displaymath}
F_{\epsilon}(\vu):= 
\left\{
\begin{array}{cl}
\displaystyle
\int_{\Omega} \bigl( \frac{\epsilon}{2} | D \vu |^2 + \frac{1}{2 \epsilon} ( 1- | \vu |^2) \bigr) dx &, \mbox{ if } \vu \in \mathcal K; \\
\infty &, \mbox{ otherwise}.
\end{array}
\right.
\end{displaymath}
Then $F_{\epsilon} \overset{\Gamma}{\rightarrow} F$, where $F$ is defined in (\ref{defF}).
\end{Theorem}
{\it Proof.} Let us first observe that for $\vu \in \mathcal K$
\begin{displaymath}
F_{\epsilon}(\vu) = \sum_{i=1}^r \int_{\Omega} \bigl( \frac{\epsilon}{2} | \nabla u_i |^2 + \frac{1}{2 \epsilon} ( u_i - u_i^2) \bigr) dx
= \sum_{i=1}^r \tilde{F}_{\epsilon}(u_i),
\end{displaymath}
where $\tilde{F}_{\epsilon}:\tilde{X}:= \lbrace v \in L^1(\Omega) \, | \, 0 \leq v(x) \leq 1 \mbox{ a.e. in } \Omega \rbrace
\rightarrow \mathbb{R} \cup \lbrace \infty \rbrace$ is defined by
\begin{displaymath}
\tilde{F}_{\epsilon}(v):= 
\left\{
\begin{array}{cl}
\displaystyle
\int_{\Omega} \bigl( \frac{\epsilon}{2} | \nabla v |^2 + \frac{1}{2 \epsilon} ( v -v^2) \bigr) dx &, \mbox{ if } v \in H^1(\Omega) \cap \tilde{X}; \\
\infty &, \mbox{ otherwise}.
\end{array}
\right.
\end{displaymath}
It is well--known (\cite{Mod87}, \cite{Alb96})  that $\tilde{F}_{\epsilon} \overset{\Gamma}{\rightarrow} \tilde{F}$ with
\begin{displaymath}
\tilde{F}(v) = \left\{
\begin{array}{cl}
\displaystyle  \frac{\pi}{8} \mathcal H^{d-1}(\partial^*\lbrace v=1 \rbrace \cap \Omega) &, \mbox{ if } 
v \in BV(\Omega, \lbrace 0,1 \rbrace); \\[2mm]
\infty &, \mbox{ otherwise}.
\end{array}
\right.
\end{displaymath}
See \cite{BloEll91,BloEll93-a,BelPaoVer90} and the following development for the calculations leading to the factor $\pi/8$.
Let $\vu \in X$ and $(\vu_{\epsilon_k})_{k \in \mathbb{N}} \subset X$ an arbitrary sequence with 
$\lim_{k \rightarrow \infty} \epsilon_k=0$ and
$\vu_{\epsilon_k} \rightarrow \vu$ in $L^1(\Omega,\mathbb{R}^r)$. Then $(u_{\epsilon_k,i})_{ \in \mathbb{N}} \subset
\tilde{X}$ and $u_{\epsilon_k,i} \rightarrow u_i$ in $L^1(\Omega), i=1,\ldots,r$, so that
\begin{displaymath}
\liminf_{k \rightarrow \infty} F_{\epsilon_k}(\vu_{\epsilon_k})  =  \liminf_{k \rightarrow \infty} \sum_{i=1}^r \tilde{F}_{\epsilon_k}
(u_{\epsilon_k,i}) \geq \sum_{i=1}^r \liminf_{k \rightarrow \infty} \tilde{F}_{\epsilon_k}(u_{\epsilon_k,i}) 
 \geq  \sum_{i=1}^r \tilde{F}(u_i) = F(\vu)
\end{displaymath}
since $\tilde{F}_{\epsilon} \overset{\Gamma}{\rightarrow} \tilde{F}$. 
It remains to show that for every $\vu \in BV(\Omega,\lbrace e_1,\ldots,e_r \rbrace) \cap X$ there exists a
sequence $(\vu_{\epsilon_k})_{k \in \mathbb{N}} \subset \mathcal K$ with $\lim_{k \rightarrow \infty} \epsilon_k=0$ such that
$\vu_{\epsilon_k} \rightarrow \vu$ in $L^1(\Omega,\mathbb{R}^r)$ and
\begin{equation} \label{limsup}
\displaystyle \limsup_{k \rightarrow \infty} F_{\epsilon_k}(\vu_{\epsilon_k}) \leq F(\vu).
\end{equation}
We essentially follow the argument in \cite{Bal90}. Because of our particular choice of potential
and the absence of volume constraints, the construction can be made more explicit allowing us at the same time to
incorporate the condition that $\sum_{i=1}^r u_i(x)=1$ a.e. in $\Omega$, which isn't considered in \cite{Bal90}. \\
Let $\vu \in BV(\Omega,\lbrace e_1,\ldots,e_r \rbrace) \cap X$, say
$\vu = \sum_{i=1}^r \chi_{E_i} e_i$. In view of Lemma 3.1 in \cite{Bal90} we can assume without loss of generality
that the $E_i$ are closed polygonal sets satisfying $\mathcal H^{d-1}(\partial E_i \cap \partial \Omega)=0, i=1,\ldots,r$.
Lemma 3.3 in \cite{Bal90} implies that there exists $\eta>0$ such that the functions $h_i: \mathbb{R}^d \rightarrow \mathbb{R}$,
\begin{displaymath}
h_i(x):= \left\{
\begin{array}{cl}
\mbox{dist}(x,\partial E_i), & x \in \mathbb{R}^d \setminus E_i, \\
- \mbox{dist}(x,\partial E_i), & x \in E_i,
\end{array}
\right.
\end{displaymath}
are Lipschitz--continuous on $H^i_{\eta}:=\lbrace x \in \mathbb{R}^d \, | \, | h_i(x) | < \eta \rbrace$ with
$| \nabla h_i(x) | =1$ a.e. in $H^i_{\eta}$. Let us introduce the function $\varphi_{\epsilon} \in C^1(\mathbb{R})$,
\begin{displaymath}
\varphi_{\epsilon}(\tau):= 
\left\{
\begin{array}{cl}
0, & \tau \leq 0; \\
\displaystyle \frac{1}{2} \bigl( 1 + \sin \bigl( \frac{\tau}{\epsilon} - \frac{\pi}{2} \bigr) \bigr), & 
0 < \tau < \epsilon \pi; \\
1, & \tau \geq \epsilon \pi.
\end{array}
\right.
\end{displaymath}
Furthermore, we define $\chi_{\epsilon}: \mathbb{R}^{r-1} \rightarrow \mathbb{R}^r$ by
\begin{displaymath}
[\chi_{\epsilon}(t)]_i:= 
\left\{
\begin{array}{cl}
1 - \varphi_{\epsilon}(t_1) &, i=1; \\
\varphi_{\epsilon}(t_1) \cdots \varphi_{\epsilon}(t_{i-1}) (1 - \varphi_{\epsilon}(t_i)) &, 2 \leq i \leq r-1; \\
\varphi_{\epsilon}(t_1) \cdots \varphi_{\epsilon}(t_{r-1}) &, i=r,
\end{array}
\right.
\end{displaymath}
where $t=(t_1,\ldots,t_{r-1})$. It is not difficult to verify that
\begin{eqnarray} 
\chi_{\epsilon}(t) &=&
\left\{
\begin{array}{ll}
e_1 &, \mbox{ if } t_1 \leq 0; \\
e_i &, \mbox{ if } t_1 \geq \epsilon \pi, \ldots, t_{i-1} \geq \epsilon \pi, t_i \leq 0;
i=2,\ldots, r-1; \\
e_r &, \mbox{ if } t_1 \geq \epsilon \pi, \ldots, t_{r-1} \geq  \epsilon \pi;
\end{array}
\right. 
\label{prop1} \\[2mm]
0 & \leq &  [\chi_{\epsilon}(t)]_i \leq 1, \, i =1,\ldots,r \; \, | D \chi_{\epsilon}(t) | \leq \frac{c}{\epsilon} \mbox{ a.e. in }
\mathbb{R}^{r-1}; \label{prop2} \\
\chi_{\epsilon}(t)& =& \frac{1}{2} \bigl( 1 - \sin \bigl( \frac{t_i}{\epsilon} - \frac{\pi}{2} \bigr) \bigr) 
e_i + \frac{1}{2}\bigl( 1 + \sin \bigl( \frac{t_i}{\epsilon} - \frac{\pi}{2} \bigr) \bigr) e_j, \label{prop3} \\[2mm]
& & \mbox{ if } 0 \leq t_i \leq \epsilon \pi, t_j \leq 0, t_k \geq \epsilon \pi, k=1,\ldots,r-1, k \neq i,j \mbox{ and } i<j.
\nonumber
\end{eqnarray}
The above function is a particular example of the function $\chi_{\epsilon}$ constructed in Lemma 3.2 in \cite{Bal90}. In addition
we have
\begin{displaymath}
\sum_{i=1}^r [\chi_{\epsilon}(t)]_i =1, \quad t \in \mathbb{R}^{r-1}.
\end{displaymath}
As a consequence, the function
$\vu_{\epsilon}(x):= \chi_{\epsilon}(h_1(x),\ldots,h_{r-1}(x)),   x \in \Omega$
belongs to $\mathcal K$ and satisfies (see p. 79 in \cite{Bal90})
\begin{displaymath}
\vu_{\epsilon} \rightarrow \vu \quad \mbox{ in } L^1(\Omega,\mathbb{R}^r), \epsilon \rightarrow 0.
\end{displaymath}
In order to analyze $F_{\epsilon}(\vu_{\epsilon})$ we introduce as in \cite{Bal90} for $i,j=1,\ldots,r$ 
the sets $\Omega^{\epsilon}_1:= E_1$,
\begin{eqnarray*}
\Omega^{\epsilon}_i &:= & \lbrace x \in E_i \, | \, h_j(x) > \epsilon \pi, j = 1,\ldots,i-1 \rbrace, \; i=2,\ldots,r; \\[2mm]
\Omega^{\epsilon}_{ij} &:= & \lbrace x \in \Omega \, | \, 0 < h_i(x) < \epsilon \pi, 
h_j(x) < 0, h_k(x) > \epsilon \pi, k \neq i,j \rbrace \mbox{ if } i<j; \\[2mm]
K^{\epsilon}_{ij} & := & \lbrace x \in \Omega \, | \,  0 \leq h_i(x) \leq \epsilon \pi,
0 \leq h_j(x) \leq \epsilon \pi \rbrace \mbox{ if } i<j.
\end{eqnarray*}
Then,
\begin{equation} \label{kij}
\displaystyle
\Omega \setminus \Bigl( \bigcup_{i=1}^r \Omega^{\epsilon}_i \cup  \bigcup_{i <j} \Omega^{\epsilon}_{ij} \Bigr) \subset
\bigcup_{i<j} K^{\epsilon}_{ij}
\end{equation}
and
\begin{equation}  \label{ueps}
\displaystyle 
\vu_{\epsilon}(x) = \left\{
\begin{array}{cl}
e_i, & x \in \Omega^{\epsilon}_i; \\[2mm]
\displaystyle \frac{1}{2} \bigl( 1 - \sin \bigl( \frac{h_i(x)}{\epsilon} - \frac{\pi}{2} \bigr) \bigr) e_i + \frac{1}{2} \bigl( 
1 + \sin \bigl( \frac{h_i(x)}{\epsilon} - \frac{\pi}{2} \bigr) \bigr) e_j, & x \in \Omega^{\epsilon}_{ij}, i<j.
\end{array}
\right.
\end{equation}
Abbreviating  $\displaystyle F_{\epsilon}(\vu,A):= \int_A \bigl( \frac{\epsilon}{2} | D \vu |^2 + \frac{1}{2 \epsilon} 
( 1 - | \vu |^2 )  \bigr) dx$ we have in view of (\ref{kij}) and (\ref{ueps})
\begin{displaymath}
F_{\epsilon}(\vu_{\epsilon}) \leq \sum_{i<j} F_{\epsilon}(\vu_{\epsilon},\Omega^{\epsilon}_{ij}) + 
\sum_{i<j} F_{\epsilon}(\vu_{\epsilon},K^{\epsilon}_{ij}).
\end{displaymath}
It is shown in \cite{Bal90} that $\limsup_{\epsilon \rightarrow 0} F_{\epsilon}(\vu_{\epsilon},K^{\epsilon}_{ij})=0$ 
for $i,j=1,\ldots,r, i<j$. Furthermore, observing (\ref{ueps}) and  
$| \nabla h_i(x)| =1$ a.e. in $\Omega^{\epsilon}_{ij}$ we obtain
\begin{displaymath}
| D \vu_{\epsilon}(x) | ^2 = \frac{1}{2 \epsilon^2} \cos^2 \bigl(\frac{h_i(x)}{\epsilon} - \frac{\pi}{2} \bigr), \; \;
1 - | \vu_{\epsilon}(x) |^2 = \frac{1}{2} \cos^2 \bigl(\frac{h_i(x)}{\epsilon} - \frac{\pi}{2} \bigr), \quad
x \in \Omega^{\epsilon}_{ij},
\end{displaymath}
so that the coarea formula yields
\begin{eqnarray*}
\lefteqn{
F_{\epsilon}(\vu_{\epsilon},\Omega^{\epsilon}_{ij})  =  \frac{1}{2 \epsilon} \int_{\Omega^{\epsilon}_{ij}} 
\cos^2 \bigl(\frac{h_i(x)}{\epsilon} - \frac{\pi}{2} \bigr) dx  = \frac{1}{2 \epsilon} \int_0^{\epsilon \pi} 
\cos^2 \bigl( \frac{t}{\epsilon} - \frac{\pi}{2} \bigr) \; \mathcal H^{d-1}(\lbrace h_i=t \rbrace 
\cap E_j) dt } \\
& = & \frac{1}{2} \int_{- \frac{\pi}{2}}^{\frac{\pi}{2}} \cos^2(s) \; \mathcal H^{d-1}
( \lbrace h_i= \epsilon (s+\frac{\pi}{2}) \rbrace \cap E_j) ds
\rightarrow \frac{\pi}{4} \mathcal H^{d-1}(\partial E_i \cap \partial E_j \cap \Omega), \epsilon \rightarrow 0.
\end{eqnarray*}
Hence,
\begin{displaymath}
\limsup_{\epsilon \rightarrow 0} F_{\epsilon}(\vu_{\epsilon}) \leq \frac{\pi}{4}  \sum_{i<j} \mathcal H^{d-1}(\partial E_i \cap \partial E_j \cap \Omega)
= \frac{\pi}{8} \sum_{i=1}^r \mathcal H^{d-1}(\partial E_i \cap \Omega) = F(\vu),
\end{displaymath}
where we note that $\partial E_i \cap \partial E_j$ is counted twice in the second sum. In conclusion,
$F_{\epsilon} \overset{\Gamma}{\rightarrow} F$. \qed

\begin{Corollary} \label{compact}
Suppose that $(\vu_{\epsilon})_{\epsilon>0} \subset \mathcal K$ is a sequence such that $(F_{\epsilon}(\vu_{\epsilon}))_{\epsilon
>0}$ is bounded. Then there exists a sequence $\epsilon_k \rightarrow 0$ and $\vu \in BV(\Omega, \lbrace e_1,\ldots,e_r \rbrace) \cap X$
such that $\vu_{\epsilon_k} \rightarrow \vu$ in $L^1(\Omega,\mathbb{R}^r)$.
\end{Corollary}
{\it Proof.} Our assumption yields that  $(\tilde{F}_{\epsilon}(u_{\epsilon,i}))_{\epsilon>0}$ is bounded for $i=1,\ldots,r$. It
is well--known that this implies that there exists
a sequence $\epsilon_k \rightarrow 0$ and $u_i \in BV(\Omega,\lbrace 0,1 \rbrace)$  such that 
$u_{\epsilon_k,i} \rightarrow u_i$ in $L^1(\Omega)$ and a.e. in $\Omega, i=1,\ldots,r$. Clearly,
$\vu_{\epsilon_k} \rightarrow \vu=(u_1,\ldots,u_r)$ in  $L^1(\Omega,\mathbb{R}^r)$, while it also follows that
$\sum_{i=1}^r u_i(x) =1$ a.e. in $\Omega$ so that $\vu \in BV(\Omega,\lbrace e_1,\ldots,e_r \rbrace) \cap X$. \qed

{\bf Acknowledgements}

The third author was supported by the EPSCR grant EP/J016780/1 and the Leverhulme Trust Grant RPG-2014-149.


\end{document}